\documentclass[envcountsame,envcountsect]{svmult}

\usepackage{makeidx}         
\usepackage{graphicx}        
\usepackage{multicol}        
\usepackage[bottom]{footmisc}

\usepackage{amsmath}
\usepackage{amscd}
\usepackage{amssymb}

  \newtheorem{thm}{Theorem}
  \newtheorem{defn}{Definition}
  \newtheorem{cor}{Corollary}
  \newtheorem{lem}{Lemma}
  \newtheorem{prop}{Proposition}
  \newtheorem{ex}{Example}
  \newtheorem{rk}{Remark}

\def\pdo{\Psi{\rm DO}}

\def\Ci{C^\infty}

\def\dd{\partial}

\def\Di{D\kern -.65em /}
\def\Dii{D\kern -.45em /}
\def\di{{\dd}\kern -.55em /}
\def\dii{{\dd}\kern -.40em /}

\def\to{\rightarrow}

\def\Dd{{\mathcal D}}

\def\Nn{{\mathcal N}}

\def\Rr{{\mathcal R}}

\def \Exp{{\rm Exp}}

\def \endsquare{ $\sqcup \!\!\!\! \sqcap$ }

\def\={\cong}
\def\>{\supset}
\def\<{\subset}

\def\12{\frac{1}{2}}
\def\2{\Dd}
\def\3{\Nn}
\def\4{\Rr}
\def\6{\cup}
\def\8{\otimes}
\def\0{^{\circ}}
\def\){\hfill{\ \qed}\enddemo}

\def\A{\Aa}

\def\R{\mathbb{R}}
\def\C{\mathbb{C}}

\def\e{\varepsilon}

\def\F{\Phi}

\def\N{\NN}

\def\Q{\QQ}
\def\Si{\Sigma}

\def\Z{\ZZ}

\def\Cl{\mbox{\rm Cl}}

\def\Ell{\mbox{\rm Ell}}

\def\Dd{{\mathcal D}}

\def\Nn{{\mathcal N}}

\def\Rr{{\mathcal R}}

\def\Si{S\kern -.65em /}

\def\tr{\mbox{\rm tr\,}}

\def\dbar{d{\hskip-1pt\bar{}}\hskip1pt}

\def \E{{\! \rm \ I \!\!\!E}}
\def \C{{\! \rm \ I \!\!\!C}}
\def \Q{{\! \rm \ I \!\!\!Q}}
\def \R {{\! \rm \ I \!R}}
\def \N {{\! \rm \ I \!N}}
\def \A {{\! \rm \ I \!\! A}}
\def \M {{\! \rm \ I \! M}}
\def \Z {{\! \rm Z\! \!Z}}

\def\cutoffint{-\hskip -12pt\int}

\def \e {{\epsilon}}

\def \End {{\rm End}}

\def \Ci {{C^\infty}}

\def \tr {{\rm tr}}
\def\pdo{\psi{\rm do}}

\begin{document}

\title*{ Chern-Weil calculus extended  to a class of  infinite dimensional manifolds}
\author{Sylvie PAYCHA \inst{}
\institute{Laboratoire de Math\'ematique, Complexe des
    C\'ezeaux, 63 177 Aubi\`ere Cedex, F
\texttt{sylvie.paycha@math.univ-bpclermont.fr}}}
%
%
\maketitle
\section*{Abstract} We discuss possible extensions of the classical Chern-Weil
formalism to an infinite dimensional setup. This is based on joint work with Steven
    Rosenberg \cite{PR1,PR2}, joint work with Simon Scott \cite{PS1,PS2} and joint work
    with Jouko Mickelsson \cite{MP}.
\section*{Acknowledgements}
I am very grateful to  Jouko Mickelsson and Steven Rosenberg and for their comments on a
preliminary version of this paper. I also thank Fr\'ed\'eric Rochon for his
useful observations and Rapha\"el Ponge for drawing my attention to relevant
references. 
\section*{Introduction}
Classical Chern-Weil formalism relates geometry to topology, assigning to the
curvature of a connection,   de Rham cohomology groups of the underlying  manifold. This theory
developped in the 40's by  Shing-Shen Chern \cite{C2} and Andr\'e Weil\footnote{In an unpublished
  paper} which can be seen as a generalisation of the Chern-Gauss-Bonnet
theorem \cite{C1},  was an important step in the theory of characteristic classes. \\
Let $G$ be a lie group with Lie algebra ${\rm Lie}(G)$.  The Chern-Weil
homomorphism assigns to  an Ad$(G)$-invariant  polynomial $f$ on ${\rm
  Lie}(G)$ a de Rham cohomology class defined as follows. Let $P\to M$ be a
$G$-principal bundle equipped with a connection $\nabla$, since the curvature
$\nabla^2$ is a ${\rm
  Lie}(G)$-valued two-form on $P$,   to a
homogeneous  Ad$(G)$-invariant polynomial $f$  of degree $j$ on ${\rm
  Lie}(G)$ corresponds a $2j$-form  $f(\nabla^2)$
on $P$. This form turns out to be  closed with de Rham cohomology class independent of the
choice of connection. \\ When $G$ is a matrix group,  Ad$(G)$-invariant
monomials on ${\rm
  Lie}(G)$ can be built from
the trace on matrices  in view of the
Ad$(G)$-invariance of the trace (section 1); the invariant polynomials are
actually generated by the monomials  $X\mapsto  {\rm tr}(X^j)$  \\
When $G$ is an infinite dimensional Lie group, there is a priori a problem to
define a trace and therefore to get invariant polynomials on ${\rm
  Lie}(G)$. We are concerned here with the Fr\'echet Lie group $\Cl^{0, *}(M,
E)$ (and its subgroups) of invertible zero order classical pseudodifferential operators acting on
smooth sections of some vector bundle $E\to M$ over a closed Riemannian
manifold $M$ (section 5). Its Lie algebra is the Fr\'echet algebra  $\Cl^{0}(M,
E)$ of zero order classical pseudodifferential operators acting on
smooth sections  $E\to M$ (section 2); it carries two types of traces \cite{LP}
together with their
linear combinations (section 3), the noncommutative residue introduced by
Adler, Manin, generalised by
Guillemin \cite{G} and Wodzicki \cite{W1} 
(see \cite{K} for a survey)  and
leading symbol traces used in \cite{PR1,PR2}. An explicit example of an
infinite rank bundle with non vanishing  first Chern class is built in
\cite{RT} using the noncommutative residue on classical pseudodifferential
operators as an Ersatz for the  trace on matrices. However, generally speaking,
Chern classes built from the noncommutative residue or leading symbol traces seem too coarse to
capture non trivial cohomology classes \footnote{Interestingly, the very fact that the class
  vanishes can  be used as a starting point to define Chern-Simons classes as in
  \cite{MRT} where the authors build non trivial Wodzicki-Chern-Simons classes
  via  a transgression of Wodzicki-Chern forms.}
so that we turn to mere linear extensions to  the algebra $\Cl^{0}(M,
E)$ of the ordinary trace on smoothing operators (section 4). 
The latter  might not be traces since they are not expected to  vanish on
brackets. We refer all the same to these as regularised traces (and weighted traces later in the
text); in constrast with the noncommutative residue and leading symbol traces
which vanish smoothing  
pseudodifferential operators, regularised traces coincide with the usual
trace on smoothing operators.   The price to pay for choosing regularised traces instead of genuine traces  is that analogues of Chern-Weil invariant
polynomials do not give rise to closed forms. Implementing techniques borrowed
from the theory of
classical pseudodifferential calculus, one  measures the obstructions to the
closedness  in terms of noncommutative residues (section 8). In specific
situations such as in hamiltonian gauge theory (section 9) where we need to
build Chern classes on pseudodifferential Grassmannians, the very locality of the
noncommutative residue can provide a way to build counterterms, and thereby to
renormalise the original non closed forms in order to turn them into closed ones. 
Loop groups \cite{F} also provide an interesting geometric setup since   
obstructions  to the closedness can vanish, thus leading to closed forms. 
On infinite rank vector bundles associated with a family of Dirac operators on
even dimensional closed spin manifolds, these obstructions can be circumvented
by an appropriate choice of regularised trace involving the very
superconnection which  gives rise to the curvature.  We discuss these last two
geometric setups in section 9. \\ \\ The paper is organised as follows: 
\begin{enumerate}
\item Chern-Weil calculus in finite dimensions
\item The algebra of (zero order) classical pseudodifferential operators
\item Traces on   (zero order)  classical $\pdo s$
\item Linear extensions of the trace on smoothing operators
\item The group of invertible zero order $\pdo s$
\item A class of infinite dimensional  manifolds
\item Singular Chern-Weil forms in infinite dimensions 
\item Weighted Chern-Weil forms; discrepancies
\item Renormalised Chern-Weil forms on $\pdo$ Grassmanians
\item Regular Chern-Weil forms  in infinite dimensions.
\end{enumerate}
\section{Chern-Weil calculus in finite dimensions}
\label{sec:1}Let $E\to X$ be a vector bundle over a $d$-dimensional manifold $X$ with
structure group $G$ a subgroup of the linear group $ {\rm Gl}_d(\C)$  and let  ${\cal A}= {\rm End}(E)$ the bundle of endormorphisms
of $E$ over $B$. Let $\Omega(X, {\cal A})$ denote the
algebra of exterior forms on $X$ with values in ${\cal A}$ equipped with the
product induced from the wedge product on forms and the product in ${\cal A}$.
If  $\sigma $ is a  section of $E$ over $X$ and $\alpha\in \Omega^k(X, {\cal
  A})$ then  $\alpha(\sigma)\in \Omega^k(X)$.   
\\
 If  $\nabla$ is a
connection on $P$ then $\nabla^2$ lies in $\Omega^2(X, {\cal A})$.  More
generally,  if $
{\cal C}(P)$ is the space of connections on $P$, to an
analytic map $f(z)$ we assign a map
\begin{eqnarray*}
f: {\cal C}(P)&\to & \Omega(X, {\cal A})\\
 \nabla &\mapsto& f(\nabla^2)=\sum_{i=0}^\infty\frac{ f^{(i)}(0)}{i!}\, \nabla^{2i}.
\end{eqnarray*}
\begin{rk}This sum is  actually finite since $\nabla^{2i}=0\quad \forall i>
  \frac{d}{2}.$
\end{rk} 
The connection $\nabla$ extends to a map
\begin{eqnarray*}
\Ci(X,TX)\times \Omega\left(X,{\cal A}\right)&\to & 
 \Omega\left(X,{\cal A}\right)\\
(U, \alpha) &\mapsto &\left(\sigma\mapsto  [\nabla_U, \alpha](\sigma):= \nabla_U (\alpha (\sigma))
+(-1)^{\vert \alpha\vert+1}\alpha(\nabla_U\sigma)\right).
\end{eqnarray*}
Here $\sigma$ stands for a section of $E$ over $X$ and $\vert \alpha\vert$ for
the degree of the form. \\
The trace ${\rm tr}: {\rm gl}_d(\C)\to \C$ on the algebra ${\rm gl}_d(\C)$   of $d$ by $d$ matrices with complex
coefficients  extends to a trace on ${\rm End}(E)$ by
\begin{eqnarray*}
{\rm tr}: {\rm End}(E)&\to & X\times \C\\
 (x,A) &\mapsto & (x,{\rm tr}(A))
\end{eqnarray*}
where ${\rm tr}$ on the r.h.s is the ordinary trace on matrices. This
is a bundle morphism  since \begin{equation}\label{eq:covarianttrace}{\rm tr}(C^{-1}\,A\,
C)={\rm tr}(A)\quad \forall C\in {\rm Gl}_d(\C), \forall A\in {\rm
  gl}_d(\C).\end{equation} \\ Similarly, to a form $ \alpha(x)= A(x)\, dx_1\wedge \cdots\wedge dx_d$ in   $\Omega\left(X,{\cal A}\right)$ corresponds a
form   ${\rm tr}(\alpha)(x):={\rm tr}(A(x))\, dx_1\wedge \cdots\wedge dx_d$ in
$\Omega(X)$.\\
From  the fact that the trace tr obeys the following properties 
\begin{equation}\label{eq:dtr} \left[d,  {\rm
  tr}\right](\alpha):=d\,{\rm tr}(\alpha)- {\rm tr} (d\,
\alpha)=0\quad  \forall \alpha\in \Omega\left(X, {\cal A}\right)\end{equation}
and 
\begin{equation}\label{eq:coboundarytr} 
\partial {\rm tr}(\alpha, \beta):= {\rm tr}\left(\alpha \wedge
  \beta+(-1)^{\vert \alpha\vert \, \vert \beta\vert} \beta \wedge
  \alpha\right)=0\quad \forall \alpha, \beta\in \Omega\left(X, {\cal A}\right),
\end{equation}  we infer the subsequent useful lemma.
\begin{lem}\label{lem:nablatr}For any $\alpha\in \Omega\left(X,{\cal
      A}\right)$
\begin{equation}\label{eq:nablatr}[\nabla, {\rm tr}](\alpha):= d\,  {\rm
    tr}(\alpha)-{\rm tr}([\nabla,\alpha])=0.
\end{equation}
\end{lem}
{\bf Proof:} In a local chart above an open subset $U$ of $X$, $$[\nabla,\alpha]= d\alpha+\theta \wedge
\alpha+(-1)^{\vert \alpha\vert+1} \alpha \wedge \theta$$ for some one
form $\theta\in \Omega^1(U, {\cal A})$ so that we can write
\begin{eqnarray*}
[\nabla, {\rm tr}](\alpha)&=& d\, {\rm tr}(\alpha)- {\rm tr}([\nabla,\alpha])\\
&=&  d\, {\rm tr}(\alpha)- {\rm tr}\left(d\alpha+\theta \wedge
\alpha+(-1)^{\vert \alpha\vert+1} \alpha \wedge \theta\right) \\
&=&- {\rm tr}\left(\theta \wedge
\alpha+(-1)^{\vert \alpha\vert+1} \alpha \wedge \theta\right)\quad {\rm
by}\quad 
{\rm  (\ref{eq:dtr})}  \\
&=& 0\quad {\rm
by}\quad 
{\rm  (\ref{eq:coboundarytr})}.
\end{eqnarray*}
\endsquare\\ \\
Combining this lemma with the Bianchi identity 
\begin{equation}\label{eq:Bianchi} [\nabla, \nabla^2]=0.\end{equation}leads to
closed Chern-Weil forms.
\begin{prop}\label{prop:ChernWeilfinitedim} For any analytic function $f$, the
  form $ {\rm tr}\left(f(\nabla^2)\right)$ is closed with de Rham cohomology class independent of the choice of
connection.
\end{prop}
{\bf Proof:} It is sufficient to carry out the proof for monomials $f(x)=x^i$
in which case we have: 
\begin{eqnarray*}
d\,{\rm tr}\left( f(\nabla^2)\right)&=&[\nabla, {\rm tr}]\left(
  f(\nabla^2)\right)+ {\rm tr}\left(
 [\nabla,  f(\nabla^2)]\right) \\
&=& {\rm tr}\left([ \nabla, \nabla^{2i}]\right)\quad {\rm by} \quad (\ref{eq:nablatr})\\
&=& \sum_{j=0}^{i}{\rm tr}\left([ \nabla, \nabla^{2}] \,
  \nabla^{2(i-1)}\right)\\
&=&0\quad {\rm by} \quad (\ref{eq:Bianchi}),
\end{eqnarray*}
which proves the closedness of ${\rm tr}\left( f(\nabla^2)\right)$. \\
Let $\nabla_t, t\in \R$ be a smooth one parameter family of connections on
$E$. Its derivative w.r. to  $t$ is a one form $\dot \nabla_t=\dot\theta_t\in
\Omega^1(X, {\cal A})$. Applying  (\ref{eq:dtr}) to $X=\R$ yields
$$[\frac{d}{dt},\,{\rm tr}]=\frac{d}{dt}\circ {\rm tr}- {\tr}\circ \frac{d}{dt}=0$$
and hence 
\begin{eqnarray*}
\frac{d}{dt}\left({\rm tr}\left( f(\nabla_t^2)\right)\right)&=& {\rm
  tr}\left( \frac{d}{dt}\nabla_t^{2i}\right) \\
&=& \sum_{j=0}^{i}{\rm tr}\left( \frac{d}{dt}\nabla_t^{2} \,
  \nabla_t^{2(i-1)}\right)\\
&=& \sum_{j=0}^{i}{\rm tr}\left(\left(\nabla_t \, \dot \nabla_t+ \dot
    \nabla_t\, \nabla_t\right)\,
  \nabla_t^{2(i-1)}\right)\\
&=& \sum_{j=0}^{i}{\rm tr}\left(\left[\nabla_t , \dot \nabla_t\right]\,
  \nabla_t^{2(i-1)}\right)\\
&=& \sum_{j=0}^{i}{\rm tr}\left(\left[\nabla_t , \dot \nabla_t\,
  \nabla_t^{2(i-1)}\right]\right)\quad {\rm by} \quad (\ref{eq:Bianchi})\\
&=&d\,  \sum_{j=0}^{i}{\rm tr}\left( \dot \nabla_t\,
  \nabla_t^{2(i-1)}\right)\quad {\rm by} \quad (\ref{eq:nablatr}).\\
\end{eqnarray*}
The variation $\frac{d}{dt}\left({\rm tr}\left( f(\nabla_t^2)\right)\right)$
is therefore exact and the de Rham class of ${\rm tr}\left(
    f(\nabla_t^2)\right)$ is independent of the parameter $t$. 
\endsquare
\section{The algebra of (zero order) classical pseudodifferential operators}
In the infinite dimensional situations considered in these notes, the algebra
of matrices ${\rm gl}_d(\C)$ on which lives the trace
used for ordinary Chern-Weil calculus, is replaced by the algebra of zero
order  classical
pseudodifferential operators on a closed manifold $M$ with values in
$\C^n$. Such an algebra   contains the algebra Map$(M, {\rm gl}_d(\C))$ of
smooth maps from $M$  to the algebra of matrices. One can think of ${\rm
  gl}_d(\C)$ as what remains of the infinite dimensional  algebra  of zero
order classical
pseudodifferential operators on $M$ when $M$ is reduced to a point
$\{*\}$. \\
We briefly  recall the definition of classical
pseudodifferential operators
($ \pdo s $) on closed manifolds, referring the reader to \cite{H},
\cite{Sh}, \cite{T}, \cite{Tr} for further details.\\
Let U be an open subset of $\R^n$. Given $a \in \C$, we consider the space of symbols
$S^a (U)$ which  consists of smooth functions $\sigma(x, \xi)$ on 
$U\times\R^n$ such that for any compact subset $K$ of
$U$ and any two multiindices $\alpha=(\alpha_1,\cdots,\alpha_n)\in \N^n, \beta=(\beta_1,\cdots,\beta_n) \in \N^n$ there exists a
constant $C_{K \alpha \beta}$ satisfying for all $(x,\xi)\in
K\times\R^n$

$$\vert \partial_x^{\alpha}\partial_{\xi}^{\beta}\sigma(x,\xi) \vert \leq
 C_{K \alpha \beta}(1+ \vert \xi \vert )^{ {\rm Re}(a)- \vert \beta\vert}.$$
where $\vert \beta\vert = \beta_1+\cdots+\beta_n.$ \\
If ${\rm Re}(a_1)<{\rm Re}(a_2)$, then $S^{a_1}(U)\subset S^{a_2}(U)$.
\begin{rk} If $a\in\R$, $a$ corresponds to the order of $\sigma\in S^a(U)$.
  The notion of order extends to complex values for  classical
  pseudodifferential symbols (see below). 
\end{rk}
The product $\star$ on symbols is defined as follows:
 if $\sigma_1 \in S^{a_1}(U)$ and $\sigma_2
\in S^{a_2}(U)$,
 $$\sigma_1\star \sigma_2(x, \xi) \sim \sum_{\alpha \in
\N^n}\frac{(-i)^{\vert \alpha\vert}}{\alpha!}\partial_{\xi}^{\alpha}\sigma_1(x,\xi)\partial_x^{\alpha}\sigma_2(x,\xi)$$
i.e. for any integer $N\geq 1$ we have
$$\sigma_1\star\sigma_2(x,\xi)
-\sum_{\vert \alpha\vert<N}\frac{(-i)^{\vert \alpha\vert}}{\alpha!}\partial_{\xi}^{\alpha}\sigma_1(x,\xi)\partial_x^{\alpha}\sigma_2(x,\xi)
\in S^{a_1+a_2-N}(U).$$
In particular, $\sigma_1\star \sigma_2 \in S^{a_1 + a_2}(U).$

 We denote by $S^{-\infty}(U):= \bigcap_{a\in \C}S^a(U)$  the algebra of
smoothing symbols on $U$ and  let  $S(U)$ be the algebra generated by $ \bigcup_{a \in \C}S^a(U)$. \\
A symbol $ \sigma\in S^a(U)$ is called  classical of order $a\in \C$  if 
$$\forall N\in \N, \quad \sigma-\sum_{j<N}\psi(\xi)\sigma_{a-j}(x, \xi) \in
S^{a-N}(U),$$
where $\sigma_{a-j}(x, \xi)$ is a positively homogeneous
 function on $U\times\R^n$ of degree $a-j$, i.e. $\sigma_{
a-j}(x,
 t\xi)=t^{a-j}\sigma_{a-j}(x, \xi)$ for all $t\in\R^+$.\\
We write for short 
\begin{equation}\label{eq:asymptsymb}\sigma(x, \xi) \sim  \sum_{j=0}^{\infty}\psi(\xi)\, \sigma_{a-j}(x, \xi).\end{equation}
Here $\psi\in C^{\infty}(\R^n)$ is any cut-off function which
vanishes for
 $\vert \xi \vert \leq \frac{1}{ 2}$ and such that $\psi(\xi)=1$ for $\vert \xi
 \vert \geq 1 .$\\   We call $a$ the order of the classical symbol $\sigma$
 and denote by  $CS^a(U)$ the subset of classical symbols of order
 $a$. The positively homogeneous component  $\sigma_a(x, \xi)$    of degree
 $a$ corresponds to the leading symbol of $\sigma$. 
\begin{ex} A smooth function $h\in \Ci(U)$ can be viewed as a
  multiplication operator $f \mapsto h\, f$ on smooth functions $f\in \Ci(U)$ and hence as a zero order
  classical symbol. 
\end{ex}
The symbol product of two classical symbols is a classical symbol
and we denote by
$$CS(U)=\langle \bigcup_{a\in \C} CS^a(U)\rangle$$ the algebra
generated by all  classical symbols on $U$.\\\
Given a symbol $\sigma\in S(U)$, we can associate to it the
continuous operator $Op(\sigma):C^{\infty}_c (U)\to C^{\infty}(U)$
defined  for $  u\in C^{\infty}_c (U)$-- the space of smooth
compactly supported functions on $U$-- by
$$\left(Op(\sigma)u\right)(x)=\int{e^{ix.\xi}\sigma(x, \xi)\widehat
  u(\xi)\dbar\xi}, $$
where   $\dbar \xi:= \frac{1}{(2\pi)^n}\, d\, \xi$ with $d\xi$ the ordinary
Lebesgue measure on $T_x^*M\simeq \R^n$ and where $\widehat u(\xi)$
 is the Fourier transform of $u$.
 Since
 $$(Op(\sigma)u)(x)=\int{\int{e^{i(x-y).\xi}\sigma(x, \xi) u(y)\dbar\xi dy}},$$
 $Op(\sigma)$ is an operator with Schwartz kernel given by
$k(x,y)=\int{e^{i(x-y).\xi}\sigma(x, \xi)\dbar\xi},$
which  is smooth off the diagonal.\\
 A pseudodifferential operator $A$  on $U$ is an operator which can be written in the form $A=Op(\sigma)+R$ where $\sigma \in S(U)$  and $R$
is a smoothing operator  i.e. $R$ has a smooth kernel.
If $\sigma$ is a classical symbol of order $a$, then $A$ is called a
classical pseudodifferential operator ($\pdo$)  of
order $a$. \\
The symbol $\sigma(A)$ of a pseudodifferential operator $A$ of order $a$ is only locally
defined whereas the leading symbol $\sigma_L(A)= \sigma_a(A)$ is globally defined.  
\begin{ex} Multiplication by $f\in \Ci(U)$ can be viewed as a zero order
  classical $\pdo$ on $U$. Here the leading symbol coincides with the symbol. 
\end{ex} 
The product on symbols
induces a composition $Op(\sigma_1*\sigma_2)=Op(\sigma_1)Op(\sigma_2).$
This in turn induces a composition on
properly supported operators.
A $\pdo$  $A$ on $U$ is called properly supported if for
any compact $C\subset U$, the set
$\{(x,y)\in {\rm Supp}(K_A), \quad x\in C  \quad {\rm or}\quad  y\in C\}$ is compact,
 where ${\rm Supp}(K_A)$ denotes the support of the Schwartz
kernel of $A$ i.e. a distribution on $U\times U$ such that, for $u
\in C^{\infty}_c (U), Au(x)=\int {K_A(x,y)u(y)dy}.$ A properly
supported $\pdo$ maps $C^{\infty}_c (U)$ into itself and admits a
symbol given by $\sigma(A)(x, \xi)=e^{-ix.\xi}Ae^{ix.\xi}$. The
composition $AB$ of two properly supported $\pdo$'s is a properly
supported $\psi$DO and $\sigma(AB)=\sigma(A) \star \sigma(B).$   \\ \\
More generally, let $M$ be a smooth  closed  manifold  of dimension $n$
and $\pi:E\to M$  a smooth  vector bundle of rank $d$ over $M$;  an operator $P:C^{\infty}
(M,E)\to C^{\infty}(M,E)$ is a  (resp. classical)
pseudodifferential operator of order $a$ if given a local
trivializing chart $(U,\phi)$ on $M$,
 for any localization
$ P_{\nu}=\chi_{\nu}^2 P\chi_{\nu}^1:C^{\infty}_c (U, \C^d)\to C^{\infty}_c
(U,  \C^d) $
 of $P$ where $ \chi_{\nu}^i \in C^{\infty}_c (U)$, the operator
$ \phi_*(P_{\nu}):= \phi P_{\nu} \phi^{-1}$ from the space $C^{\infty}_c
(\phi(U), \C^d)$
into $C^{\infty} (\phi(U), \C^d)$ is a  (resp classical)
pseudodifferential operator of order $a$.
\begin{ex} A smooth section $f\in \Ci(M, {\rm  End}(E))$ can be viewed as a multiplication
  operator $u\mapsto f\, u$ on smooth sections  $u$ of $E$ and hence as a zero
  order classical $\pdo$
. 
\end{ex}
 Let $\Cl^a(M, E)$ denote
the set of classical pseudodifferential operators of order $a$.\\
If $A_1\in \Cl^{a_1}(M, E), A_2\in \Cl^{a_2}(M, E)$,
 then $A_1 A_2\in \Cl^{a_1 + a_2}(M, E)$
and we denote by
$$\Cl(M,E):= \langle\bigcup_{a\in \C} \Cl^a(M, E)\rangle$$
the algebra generated by all classical pseudodifferential operators acting
on smooth sections of $E$.\\ It follows from the above discussion that 
$$\Ci\left(M, {\rm End}(E)\right)\subset \Cl^0(M, E)\subset \Cl(M, E).$$
\begin{rk}  When $E$ is the trivial bundle
$M\times \R$, we  drop $E$ in the notation writing $\Cl^a(M)$,
$\Cl^{-\infty}(M)$, $\Cl(M)$
instead of  $\Cl^a(M, E)$, $\Cl^{-\infty}(M,E)$, $\Cl(M,E)$.
\end{rk}
\begin{rk}When $M$ reduces to a point, then $E$ is a vector space corresponding to the
model space of the original bundle $M$ and we have
$$M=\{*\}\Longrightarrow \Ci\left(M, {\rm End}(E)\right)= \Cl(M, E)=\Cl^0(M,
E)=\Cl^{-\infty}(M, E)= {\rm End}(E).$$
 \end{rk}

\section{Traces on   (zero order)  classical $\pdo s$}
Having chosen $\Cl^0(M, E)$ as a potential infinite dimensional  Ersatz for the algebra
${\rm gl}_d(\C)$, it remains to find  linear forms on  $\Cl^0(M, E)$  as an
Ersatz for the trace on matrices. \\
The ordinary trace on matrices extends to a trace on smoothing operators:
\begin{eqnarray*}{\rm tr}: \Cl^{-\infty}(M, E)&\to & \C\\
A&\mapsto &  \int_M {\rm tr}_x(k_A(x,x))\, dx=\int_{T^*M}{\rm
  tr}_x(\sigma(A)(x, \xi))\, dx\,\dbar
\xi,
\end{eqnarray*}
where  $k_A$ stands for the Schwartz  kernel of $A$, $\sigma_A$
for the symbol of $A$ and ${\rm tr}_x$ for the fibrewise trace defined previously
using the ordinary trace on matrices.
\\ But it does not further  extend to a trace on $\Cl^0(M, E)$  i.e. to  a linear
form $\lambda: \Cl^0(M, E)\to \C$ which vanishes on brackets $$\partial
\lambda(A, B):= \lambda\left([A, B]\right)=0\quad \forall A, B\in \Cl^0(M, E).$$
A  well known result by Wodzicki  \cite{W1} (see also \cite{W2} and \cite{K}
for a review) and
proved independently by Guillemin \cite{G} gives the uniqueness (up to a
multiplicative factor) of a trace on the whole algebra 
$\Cl(M, E)$ of classical pseudodifferential operators.\footnote{  Since then
  other proofs, in particular a homological proof on symbols
in \cite{BG}  (see also
\cite{P2} for another alternative proof)  and various extensions
of this uniqueness result were derived, see \cite{FGLS} for a
generalisation to manifolds with boundary, see \cite{S} for a generalisation to
manifolds with conical singularities (both of which prove uniqueness up to
smoothing operators),  see \cite{L} for an extension to
log-polyhomogeneous operators as well as for an argument due to Wodzicki to
 get uniqueness on the whole algebra of classical operators, see \cite{Po2}
 for an extension to Heisenberg manifolds.}  \\
Indeed, Wodzicki showed that any trace on $\Cl(M, E)$ is proportional to the noncommutative residue defined as follows.  The residue density at point $x\in M$
$$\omega_{\rm res}(A)(x):= \left(\int_{S_x^* M} {\rm tr}_x\left(\sigma_{-n}(A)(x,
  \xi)\right)\, \dbar_S \xi\right)\, dx$$
where $S_x^*M \subset T_x^*M$ is the cotangent unit sphere (here $\dbar_S\xi:=\frac{d_S\xi}{(2\pi)^n}$
  stands for the (normalised) volume measure on the cotangent unit sphere $S^*M$
  induced by the canonical volume measure on the cotangent bundle
  $T^*M$ and $(\cdot)_{-n}$ denotes the
positively homogeneous component  degree $-n$ of the symbol) is globally
defined so that the {\bf noncommutative residue} \footnote{It generalises to
higher dimensions a notion of residue previously introduced by Adler and Manin in the one
dimensional case.}
\begin{equation}\label{eq:Wodzickiresidue}{\rm res}(A):= \int_M \omega_{\rm
    res}(A)(x):= \int_M dx\, {\rm res}_x(A)
\end{equation}
is well defined on $\Cl(M, E)$. 
 \\ \\ Restricting to zero order classical pseudodifferential operators allows for
another type of trace, {\bf leading symbol traces}
associated with any linear form  $\tau$ on $ C^\infty(S^*M)$  introduced in
\cite{PR1,PR2} in relation to Chern-Weil forms:
  $${\rm Tr}_0^{\tau}(A):= \tau \left({\rm
        tr}_x\sigma_0(A)(x,\xi)\right).$$
 Whenever $\tau(1)\neq 0$ we set:
\begin{equation}\label{eq:leadingsymboltrace}{\rm tr}_0^{\tau}(A):=\frac{{\rm Tr}_0^{\tau}(A)}{\tau \left(1\right)}=\frac{\tau \left({\rm
        tr}_x\left(\sigma_0(A)(x,\xi)\right)\right)}{ \tau \left(1\right) }.
\end{equation}
\begin{thm}\label{thm:LP}\cite{LP}
 All traces on the algebra $\Cl^{ 0}(M,E)$ are linear
  combinations the Wodzicki residue and leading symbol traces. 
\end{thm} 
\begin{rk}When $M$ reduces to a point so that $n=0$ and $A$ is a matrix,  
then both  res$(A)$ and ${\rm tr}_0^\tau(A)$  are
proportional to  the ordinary matrix trace.
\end{rk}
 Both   the noncommutative
residue and the leading symbol traces clearly vanish  on smoothing operators
and therefore neither of them   extends the
ordinary trace on smoothing operators. 
If we insist on building linear forms on $\Cl^0(M, E)$ that extend the ordinary
trace on smoothing operators, we need to drop the requirement that it vanishes
on brackets. The linear forms we are about to describe are actually defined on
the whole algebra $\Cl(M, E)$.\\
We use the unique extension \cite{MSS} (see also \cite{P2} where the
uniqueness of the noncommutative residue and the canonical trace are handled
simultaneously), called the canonical trace, of the trace on smoothing operators to
the set $$\Cl^{\notin \Z} (M,E):= \bigcup_{a\notin \Z}\Cl^a(M,E)$$ of {\it non integer order operators} in
$\Cl(M, E)$. It was  popularised by Kontsevich and Vishik in \cite{KV} even though it
was known long before by Wodzicki and Guillemin and is defined as follows.
\\ For any $A\in \Cl(M, E)$, for
any $x\in M$, one can infer from (\ref{eq:asymptsymb}) (see e.g. \cite{L}) that the integral $\int_{B_x(0, R)} {\rm
  tr}_x\sigma_A(x, \xi)\dbar \xi $ of the fibrewise trace ${\rm tr}_x
\sigma_A$ of its symbol
$\sigma_A$  over the ball $B_x(0,R )$
of radius $R$ and centered at $0$ in the cotangent bundle $ T_x^*M$, has an  asymptotic expansion in
decreasing powers of $R$ which is polynomial in  log $R$. Picking the constant
term yields the cut-off
regularised integral $$\cutoffint_{T_x^*M}{\rm
  tr}_x\left(\sigma(A)(x, \xi) \right)\,\dbar \xi := {\rm fp}_{R\to \infty}\int_{B_x(0, R)} {\rm
  tr}_x\left(\sigma(A)(x, \xi)\right)\,\dbar \xi$$ which clearly  coincides with the
ordinary integral on smoothing symbols. 
\begin{thm} \cite{KV} Whenever the  operator $A\in \Cl(M, E)$ has non integer order or has order $<-n$
then
$$\omega_{KV}(A)(x):= \left(\cutoffint_{T_x^*M} {\rm tr}_x\left(\sigma(A)(x,
  \xi)\right)\, \dbar \xi\right)\, dx \quad \forall x\in M,$$
defines a global density on $M$ so that the  canonical trace
\cite{KV} (see also \cite{L} for an extension to log-polyhomogeneous operators):
$${\rm TR}(A):= \int_M \omega_{KV}(A)(x):= \int_M {\rm TR}_x(A) \, dx$$
makes sense \footnote{ However, in  general $\omega_{KV}(A)(x)$ is only locally defined
  and does not integrate over $M$ to a well defined linear form.}.
 The canonical trace  vanishes on brackets of non integer
order or of order $<-n$ \cite{KV} (see also \cite{L}) i.e.
$${\rm TR}\left([A, B ]\right)= 0\quad \forall A, B\in \Cl(M, E)\quad {\rm
    s.t.}\quad [A, B]\in \Cl^{\notin \Z}(M, E)\cup \Cl^{<-n}(M, E).$$
\end{thm}
\begin{rk}
For any smoothing operator $${\rm TR}(A)=\cutoffint_{T_x^*M} {\rm tr}_x\left(\sigma(A)(x,  \xi)\right)\, \dbar \xi\, dx ={\rm tr}(A)$$ so that the canonical
trace indeed extends the ordinary trace on smoothing operators. 
\end{rk}

\section{Linear extensions of the trace on smoothing operators}
Unfortunately, the operators one comes across in infinite dimensional geometry
as well as in quantum field theory are typically integer order operators such
as the Laplace operator, the Dirac operator, the Green operator...so that we
cannot implement the canonical trace on such operators.\\
In order to match the canonical trace with our needs in spite of this apparent 
discrepancy, we perturb the operators holomorphically $A\mapsto A(z)$ thereby perturbing their
order $a\mapsto \alpha(z)$  and we define regularised trace of such operators
as finite parts at $z=0$ of TR$(A(z))$. \\ To carry out this construction we need the notion of holomorphic
family of symbols which we now recall. 
\begin{defn}
 Let $\Omega$ be a domain of $\C$. A family $(\sigma(z))_{z\in \Omega}\subset CS(U)$ is holomorphic when \\
(i) the order $\alpha(z)$ of $\sigma(z)$ is holomorphic  on $\Omega$.\\
(ii) For $(x, \xi)\in U\times\R^n$, the function $z\to\sigma(z)(x,
\xi)$ is holomorphic on $\Omega$ and
 $\forall k \geq 0, \partial_z^k\sigma(z)\in S^{\alpha(z)+ \epsilon}(U)$ for any $\epsilon>0.$\\
(iii) For any integer $j\geq 0,$ the homogeneous symbol
$\sigma(z)_{\alpha(z)-j}(x, \xi)$ is holomorphic on $\Omega.$
\end{defn}
It leads to the following  notion of holomorphic family of $\pdo s$. \\
 \begin{defn} A  family $z\mapsto A(z)\in \Cl^{\alpha(z)}(M, E)$ of
   log-classical$ \pdo s$ parametrized by a domain $\Omega\subset \C$ is holomorphic
if in each local trivialisation of $E$ one has $$A(z) = {\rm
Op}(\sigma_{A(z)}) + R(z)$$ with $\sigma_{A(z)}$  a holomorphic
family of classical symbols of order $\alpha(z)$\footnote{In applications the
  order is affine in $z$.} and $R(z)$ a smoothing
operator with Schwartz kernel $R(z,x,y)\in \Ci(\Omega\times M\times
M,\End(V))$ holomorphic in $z$ where $V$ is the model space of $E$.
\end{defn}
A holomorphic family of classical operators of holomorphic order
$\alpha(z)$ parametrised by $\Omega$ has integer order $\geq -n$  on the set  $\Omega\cap
\alpha^{-1}\left(\Z\cap [-n,\infty[\right)$.  Outside that set, the canonical
trace $ {\rm TR}\left( A(z)\right)$ is therefore well defined. 
\begin{thm}\label{thm:KVop} \cite{KV}\cite{L}
 Let $z\mapsto A(z)\in
 \Cl^{\alpha(z)}(M, E)$
be a holomorphic family of classical $\pdo$s  on a domain $\Omega\subset \C$. Then
the map
$$z\mapsto {\rm TR}\left( A(z)\right)$$
is meromorphic with   poles  of order $1$ at points $z_j\in \Omega\cap
\alpha^{-1}\left([-n, +\infty[\, \cap \, \Z\right)$ such that
$\alpha^\prime(z_j) \neq 0$. 
\end{thm}
\begin{defn} A  holomorphic {\rm regularisation  scheme }  on  $ \Cl(M, E)$ is
  a linear map which sends 
$  A\in  \Cl(M, E)$ to  a holomorphic family $A(z)\in \Cl(M, E)$ 
such that $A(0)=A$ and $A(z)$ has order $\alpha(z)$ with $\alpha$ holomorphic and
$\alpha^\prime(0)\neq 0$.
\end{defn}
In order to illustrate this with examples, it is useful to introduce the
following definition. \begin{defn}An operator $A \in \Cl(M,E)$ has principal angle $\theta$ if for every $(x,
\xi)\in T^*M-\{0\}$, the leading symbol $\sigma_A^L(x, \xi)$ has no
eigenvalues on the ray  $\overline{L_\theta}=\{re^{i\theta}, r\geq 0\}$; in that case $A$ is elliptic. \\ We call an  operator $A\in \Cl(M, E)$
 admissible with spectral
 cut $\theta$
if $A$ has principal angle $\theta$ and the spectrum of $A$ does
 not meet the open ray  $ {\text{L}}_{\theta}=\{re^{i\theta}, r>
 0\}$. In particular such an operator is elliptic, i.e. it has invertible
 leading symbol $\sigma_L(A)(x, \xi)\in{\rm End}(E_x)$ for
 all $x\in M, \xi\in T_x^*M-\{0\}$ where $E_x$ is the fibre of $E$ over $x$.
\end{defn}
\begin{rk}
 When $A$ has principal angle $\theta$ and the spectrum of $A$ does not meet
$L_{\theta}$, $\theta$ is called an Agmon angle of $A$. In that case, $A$ is
invertible elliptic. We want to allow for non invertible operators hence the
need for introducing admissibility.
\end{rk}
\begin{ex} $\zeta$- regularisation
\begin{equation}\label{eq:zetareg}{\cal R}:A\mapsto A(z):= A\, Q_\theta^{-z}\end{equation} with $Q$ an admissible operator in
$\Cl(M, E)$ with positive order $q$ and spectral cut $\theta$ yields   typical (and very useful)
examples of holomorphic regularisations.\end{ex}
On the basis of the results of the previous section, given a holomorphic
regularisation ${\cal R}: A\mapsto A(z)$,  we can pick the finite part in the Laurent
expansion ${\rm TR}\left(A(z)\right)$ and 
 set the following definition.
\begin{defn} A holomorphic regularisation scheme ${\cal R}:A\mapsto A(z) $ on
 $\Cl(M, E)$ 
induces a linear form:
\begin{eqnarray*}
{\rm tr}^{\cal R}:\Cl(M, E)& \to &\C\\
A &\mapsto & {\rm tr}^{\cal R} (A):= {\rm fp}_{z=0}
{\rm TR}\left(A(z)\right)
\end{eqnarray*} 
called ${\cal R}$-regularised trace \footnote{It carries this name because it
  extends the ordinary traace on smoothing operators and in spite of
  the fact that it does not vanish on brackets as we shall soon see.}.
When  ${\cal R}$ is a $\zeta$-regularisation (\ref{eq:zetareg}) determined by
an admissible operator $Q$ with spectral cut $\theta$  we call 
$Q$-weighted trace the linear form  ${\rm tr}^{\cal R} $ and denote it by
${\rm tr}_\theta^Q$. 
\end{defn}
The following result measures the difference between the regularised trace and
the (generally non existing) canonical trace. 
\begin{thm} Let ${\cal R}: A \mapsto A(z)$ with $A(z)$ of order $\alpha(z)$ be
  a holomorphic regularisation with order $z\mapsto \alpha(z)$ affine in $z$.\\
The linear form ${\rm tr}^{\cal R}$ extends the usual trace defined on   operators of
order $<-n$ as well as the canonical trace TR defined on non integer order
operators) to $\pdo$s of all orders.
\\ Moreover \cite{PS1}, 
\begin{equation}\label{eq:RegPS}
 {\rm tr}^{\cal R}\left(A\right)
=\int_M\,\dbar x\ \left({\rm TR}_x(A)- 
\frac{1}{\alpha^\prime(0)} {\rm res}_{ x}
\left( A^{\prime}(0)\right)\right)
\end{equation} 
where, in spite of the fact that $A^\prime(0)$ is not anymore expected to be
classical \footnote{It is log-polyhomogeneous of $\log$ type $1$ \cite{PS1}, meaning by
  this that the asymptotic expansion (\ref{eq:asymptsymb}) might present a
  logarithmic divergence $\log \vert \xi\vert$ in $\vert\xi \vert$ as $\vert
  \xi\vert \to \infty$.}, its residue density is defined in a similar manner: $${ \rm
res}_{x}( A^{\prime}(0)):= \int_{S^*M} {\rm tr}_x(\sigma(A))_{-n}(x,
\xi)\, \dbar\xi.$$
\end{thm}
\begin{rk} When the residue density ${\rm res}_{ x}
\left( A^{\prime}(0)\right)$ vanishes,   ${\rm TR}_x(A)\, dx$ defines a global
density and ${\rm tr}^{\cal R}\left(A\right)
={\rm TR}(A)$.   In particular    $$ {\rm tr}^{\cal R}\left(A\right)
={\rm TR}(A)\quad \forall A\in \Cl^{\notin \Z}(M, E)$$ is independent of the
regularisation scheme. 
\end{rk}
 $\zeta$-regularisations provide an interesting class of examples. Let 
$${\cal R}:A\mapsto A(z):= A\, Q_\theta^{-z}$$ with $Q$ an admissible operator in
$\Cl(M, E)$ with positive order $q$ and spectral cut $\theta$. Then
$A^\prime(0)= -A\, \log_\theta Q$. Here $\log_\theta Q$ stands for the logarithm  of an admissible operators $Q\in  \Cl(M, E)$ with
spectral cut $\theta$  defined in terms of the derivative at $z=0$
of its complex power \cite{Se}:
$$\log_\theta Q=  \partial_z {Q_\theta^{z}} _{\vert_{z=0}},$$ where
$Q_\theta^z$ is the complex power of $A$ defined using a Cauchy integral on a
contour $\Gamma_\theta$ around the spectrum of $A$. Formula (\ref{eq:RegPS})
therefore  reads:
\begin{equation}\label{eq:zetaPS}
 {\rm tr}_\theta^Q\left(A\right)
=\int_M\,\dbar x\ \left({\rm TR}_x(A)- 
\frac{1}{q} {\rm res}_{ x}
\left( A\, \log_\theta Q\right)\right).
\end{equation}
We borrow the following definition  from \cite{OP}.\begin{defn}\label{eq:condtraces} We call an operator $A\in \Cl(M, E)$ conditionally trace-class whenever
the fibrewise trace of
its symbol ${\rm tr}_x\left(\sigma_A(x, \cdot)\right)$ at a point $x\in M$
is of order $<-n$  in which case we set
$${\rm tr}_{\rm cond}(A):= \int_M dx\, \int_ {T_x^*
  M}{\rm tr}_x  \left(\sigma(A)(x, \xi)\right)\, \dbar\xi$$
 which we call the conditioned trace of $A$.
\end{defn}
\begin{ex} Clearly, operators in $\Cl(M, E)$ of order $<-n$  are conditionally
  trace-class and their conditioned
  trace coincides with their ordinary trace.
 \end{ex}
\begin{ex} Let $E=M\times \R^n$ be a rank $n$ trivial vector bundle over $M$ then
$({\cal A}u)(x):= A\, u(x)\quad\forall x\in M, u\in C^\infty (M,
\R^n)$ with $A\in o(\R^n)$ (the Lie algebra of the orthogonal group
$O(\R^n)$) is not trace class since it is a multiplication operator.
However, it is conditionally trace-class with zero conditioned trace
since the fibrewise trace of its symbol which coincides with the
trace of the matrix $A$, vanishes.
\end{ex}
\begin{ex} More generally,  let us consider a trivial vector
  bundle $E= M\times V$ with $V$ a finite dimensional space, then $\Cl(M, E)\simeq \Cl(M)\otimes {\rm Hom}(V)$ so
  that
  ${\rm tr}_V(A) \in \Cl(M)$ for any $A\in \Cl(M, E)$ where ${\rm tr}_V$ is
  the ordinary trace on  Hom $(V)$.  If  ${\rm tr}_V(A)$ is trace-class then $A$ is
  conditionally trace-class.
\end{ex}
\begin{prop} \label{prop:condtrace}
Let $A\in \Cl(M, E)$ be conditionally trace-class and ${\cal R}: A\mapsto
A(z)$ a holomorphic regularisation. Then
\begin{enumerate}
\item ${\rm res}(A^\prime(0))$ vanishes, 
\item $A$ has a well
  defined canonical trace  $${\rm TR}(A):= \int_M  \dbar x\,  \int_{T_x^*M} {\rm tr}_x\left(\sigma(A)(x,
  \xi)\right)\, \dbar \xi$$ and
$${\rm tr}^{\cal R}(A)= {\rm tr}_{\rm cond}(A)={\rm TR}(A).$$
\end{enumerate}
\end{prop}
{ \bf Proof:} The assertions follow from  the fact that the scalar symbol
${\rm tr}_x(\sigma_A(x,\cdot))$ together with the fact that the derivative \footnote{ We recall that
    $A^\prime(0)$ has the same
  order as $A$.}  ${\rm
    tr}_x(\sigma_{A^\prime(0)}(x,\cdot))$
  are of order $<-n$.\\
Indeed, this implies that ${\rm res}_x\left( A^\prime(0)\right)$ vanishes and
hence that ${\rm tr}^{\cal R}(A)= {\rm TR}(A)$ by (\ref{eq:RegPS}).
\endsquare
\section{The group of invertible zero order $\pdo s$}
The Lie  algebra $\Cl^0(M, E)$ offers a natural generalisation of the 
algebra ${\rm End}(E)$. The corresponding Lie group of invertible zero order
$\pdo$'s offers a   natural generalisation of the group GL$(E)$ of
linear transformations of a vector space $E$.  \\ \\  For $a\in \C$, the
linear space $\Cl^{a}(M, E)$ of classical pseudodifferential operators of order
$m$ can be equipped with a Fr\'echet topology. For this, one first
equips the set $CS^{a}(U, W)=CS^a(U)\otimes {\rm End}(W)$ of classical symbols of order $a$ on an
open subset $U$ of $\R^n$ with values in an euclidean vector space $W$
(with norm $\Vert\cdot \Vert$) with a Fr\'echet structure.    The following
semi-norms labelled by multiindices $\alpha,\beta$ and integers $j\geq
0$, $N$ give rise to a Fr\'echet topology on $CS^m(U, W)$  (see \cite{H}):
\begin{eqnarray*}
&{} & {\rm sup}_{x\in K, \xi \in \R^n} (1+\vert \xi\vert)^{-{\rm Re}(a)+\vert \beta\vert} \, \Vert \partial_x^\alpha \partial_\xi^\beta \sigma(x, \xi)\Vert;\\
&{}&  {\rm sup}_{x\in K, \xi\in \R^n}  (\vert \xi\vert)^{-{\rm Re}(a)+N+\vert \beta\vert}\Vert \partial_x^{\alpha} \partial_\xi^{\beta}\left(\sigma-\sum_{j=0}^{N-1} \psi(\xi)\, \sigma_{a-j}\right)(x, \xi) \Vert;\\
&{}& {\rm sup}_{x\in K, \vert\xi\vert=1}  \Vert \partial_x^{\alpha} \partial_\xi^{\beta} \sigma_{m-j}(x, \xi) \Vert,
\end{eqnarray*}
where $K$ is any compact set in $U$.\\
Given a vector bundle based on a closed manifold $M$, the set
$\Cl^a(M, E)$ of classical pseudodifferential operators acting on
sections of $E$, namely  pseudodifferential operators $A$ acting on
sections of $E$ that have local classical symbols $\sigma^U(A)\in CS(U, W)$ in a local trivialization $E_{\vert_U}\simeq U\times W$, inherits a
Fr\'echet structure via the Fr\'echet structure on classical symbols.
Given an atlas $(U_i, \phi_i)_{i\in I}$ on $M$ and corresponding local
trivializations $E_{\vert_{U_i}}\simeq U_i \times W$ where $W$ is the model
fibre of $E$, in any local chart $U_i$,
we can equip $\Cl^{a}(M, E)$ with the following family of seminorms
labelled by multiindices $\alpha,\beta$ and integers $j\geq 0, i\in
I$, $N\geq 0$
\begin{eqnarray*}
&{} & {\rm sup}_{x\in K, \xi \in \R^n} (1+\vert \xi\vert)^{-{\rm Re}(a)+\vert \beta\vert} \, \Vert \partial_x^\alpha \partial_\xi^\beta \left(\sigma^{U_i}(A)\right)(x, \xi)\Vert;\\
&{}&  {\rm sup}_{x\in K, \xi\in \R^n}  (\vert \xi\vert)^{-{\rm Re}(a)+N+\vert \beta\vert}\Vert \partial_x^{\alpha} \partial_\xi^{\beta}\left(\sigma^{U_i}(A)-\sum_{j=0}^{N-1} \psi(\xi)\, \sigma_{a-j}^{U_i}(A)\right)(x, \xi) \Vert;\\
&{}& {\rm sup}_{x\in K, \vert\xi\vert=1}  \Vert \partial_x^{\alpha} \partial_\xi^{\beta} \sigma_{a-j}^{U_i}(A)(x, \xi) \Vert,
\end{eqnarray*}
where $K$ is any compact subset of $\phi_i(U_i)\subset \R^n$.\\
\begin{prop} $\Cl^{ 0}(M, E)$ is  a Fr\'echet Lie algebra and the traces of
  Theorem \ref{thm:LP} are continuous for the Fr\'echet topology. 
\end{prop}
{\bf Proof:} The continuity of the traces can easily be seen from their very
definition. Let us discuss the continuity of the bracket. Since $\sigma(AB)\sim \sigma(A)\star \sigma(B)$ for two operators $A,
B$ with symbols $\sigma(A)$, $\sigma(B)$,  the product
map on $\Cl^0(M, E)$ is smooth as a consequence of the smoothness of the symbol product
$\sigma\star \tau \sim\sum_\alpha
\frac{1}{\alpha!}\partial_\xi^\alpha \sigma \, \partial_x^\alpha
\tau$ 
 on $CS(U, V)$ for any vector space $V$. It follows that the bracket  is a
 continuous bilinear map on
 $\Cl^0(M, E)$.\endsquare   
\\ \\
Let $${\Cl}^{0,*}(M, E):= \{A\in \Cl^0 (M, E), \quad
\exists A^{-1}\in \Cl^0 (M, E)\} $$
be the group of invertible zero
order classical pseudodifferential operators which is strictly contained in
the intersection  $\Cl^0(M, E)\cap \Cl^*(M, E)$
where 
$$\Cl^*(M, E)= \{A\in \Cl(M, E) ,\quad  \exists \, A^{-1}\in \Cl(M, E)\}$$
is the group of invertible classical pseudodifferential operators.
 \begin{rk} It is useful to note that $\Cl^{*}(M, E)$ acts on $\Cl^{ a}(M, E)$ for any $a\in \C $ by the adjoint action defined for $P\in {\Cl}^{*}(M, E)$ by
\begin{eqnarray}\label{eq:adjointaction}
\Cl^{a}(M,E)&\to & \Cl^{ a}(M,E)\nonumber\\
A&\mapsto& Ad_P A:= P^{-1}A P
\end{eqnarray}
and specifically  on the algebra  $\Cl^0(M, E)$.
\end{rk}

\begin{prop}\cite{KV}
  $\Cl^{0, *}(M, E)$ is a Fr\'echet Lie group  with  Lie alegbra  $\Cl^{ 0} (M, E)$.
\end{prop} 
{\bf Proof:} We only discuss the continuity of the inverse map, referring the
reader to \cite{KV} for further details.\\ 
As an open subset in the Fr\'echet space $\Cl^0(M, E)$, ${\Cl}^{0,*}(M, E)$ is a Fr\'echet manifold modelled on $\Cl^0(M, E)$.\\
We already know that the product map is smooth. The smoothness of the inversion $A\mapsto A^{-1}$
on $\Cl^{0, *}(M, E)$ follows from the fact that for an operator $A\in
\Cl^{0, *}(M, E)$ with symbol $\sigma(A)$ and order $a$, the positively
homogeneous components of its inverse $A^{-1}$ of order $-a$ are given by
$$\left(\sigma(A^{-1})\right)_{-a-j}=\frac{i}{2\pi}\int_{\Gamma}\lambda^{-1}
\, \tau_{-a-j}(A) \, d\lambda$$
where $\Gamma$ is a contour around the
spectrum of $A$ and where
 \begin{eqnarray*}
\tau_{-a}(A)&=& (\sigma(A)-\lambda)^{-1},\\
 \tau_{-a-j}(A)&=& -\tau_{-a}(A)\, \sum_{k+l+\vert \alpha\vert =j, l,j} i^{-\vert \alpha\vert} \frac{1}{\alpha!} \, \partial_\xi^\alpha \sigma_{a-k}(A)\partial_x^\alpha \tau_{-a-l}(A).
\end{eqnarray*}\endsquare\\ \\
Following \cite{KM} we say that a  Lie group ${\mathcal G}$ admits an exponential mapping if there exists
  a smooth mapping $$\Exp: {\rm Lie}\left({\mathcal G}\right)\to {\mathcal
    G}$$
  such that $t\mapsto \Exp (t\, X)$ is a  one-parameter subgroup with tangent vector $X$. $\Exp(0)=
  e_{\mathcal G}$ and $\Exp$ induces the identity map $D_e \Exp=Id_{{\rm
      Lie}({\mathcal G})}$ on the corresponding Lie algebra.\\
 All known smooth Fr\'echet Lie groups and  in particular the group $\Cl^{0,*}(M,
E)$ (see \cite{KV}) admit an exponential
  mapping  although it is not known, according to \cite{KM}, whether
  any smooth Fr\'echet Lie group does admit an exponential mapping.
\\ \\ The topology of $\Cl^{0, *}(M, E)$ has been investigated in various
contexts. \\
Recall  (see e.g. \cite{Ka}) that the
  fundamental group $\pi_1(GL_d(\C))$ is generated by the homotopy
  classes $[l]$ of the loops $$l(t)= e^{2i\pi t\, \pi}$$ where $\pi: \C^d\to
  \C^d$ is a projector. \\
 A similar statement holds for the fundamental group of
  $\Cl^{0, *}(M, E)$ with the projectors $\pi$ replaced by  pseudodifferential
operators introduced by Burak \cite{Bu}, later used by Wodzikci \cite{W1} and
further investigated by Ponge \cite{Po1}  which encode  the spectral asymmetry of elliptic classical
pseudodifferential operators:
$$\Pi_{\theta, \theta^\prime} (Q):= \frac{1}{2i\pi} \int_{C_{\theta,
    \theta^\prime}} \lambda^{-1} \, Q \, (Q-\lambda)^{-1} \,
d\lambda$$
where $$C_{\theta, \theta^\prime}:= \{ \rho \,
e^{i\theta}, \infty > \rho\geq r\} \cup \{r\, e^{i\,t}, \theta\leq
t\leq \theta^\prime\} \cup \{\rho\, e^{i\theta^\prime}, r\leq \rho
<\infty\},$$
with $Q\in \Cl(M, E)$ elliptic with positive order and whereby
$r$ is chosen small enough so that non non-zero eigenvalue of $Q$ lies in
the disc $\vert \lambda\vert \leq r$. It turns out that $\Pi_{\theta,
  \theta^\prime}(Q)$ is a bounded $\pdo$  projection  on $L^2(M, E)$ (see
\cite{BL} and
\cite{Po1})
and either a zero'th order pseudodifferential operator or a smoothing
operator. For $Q$ of order $q$ with leading symbol $\sigma_L(Q)$, the leading
symbol of $\Pi_{\theta, \theta^\prime} (Q)$ reads:$$\pi_{\theta, \theta^\prime} (\sigma_L(Q)):= \frac{1}{2i\pi}
\int_{C_{\theta, \theta^\prime}} \lambda^{-1} \, \sigma_L(Q)\,
(\sigma_L(Q)-\lambda)^{-1} \, d\lambda.$$
\\ \\ The following proposition (see \cite{KV}, see also \cite{LP}) shows that
these pseudodifferential projectors generate the fundamental group
$\pi_1\left(\Cl^{0,*}(M,E)\right)$. Let  ${\rm GL}_\infty( {\mathcal A})$ be the direct limit \footnote{A natural
  embedding ${\rm GL}_n(A)\to {\rm GL}_{n+1}(A) $ of  an $n\times n$
matrix $g\in {\rm GL}_{n}(A)$ in ${\rm GL}_{n+1}(A)$  is obtained  inserting
$g$ in the upper left corner, $1$ in the lower right corner
and filling the other slots in the  last line and column with zeroes.} of linear groups
${\rm GL}_n ({\cal A})$. \\ \\
\begin{prop}
\item $\pi_1\left({\rm GL}_\infty\left(\Cl^{0,*}(M,E)\right)\right)$ is
  generated by the homotopy class of loops
   \begin{equation}\label{eq:generators}L_{\theta,
         \theta^\prime}^Q(t):=e^{2i\pi\, t \Pi_{\theta, \theta^\prime}(Q)}
\end{equation} where $Q \in \Cl(M,
  E)$ is any elliptic operator with positive order.
\end{prop}
\begin{rk} When $M$ reduces to a point $\{*\}$ then $\sigma_L(Q)$ reduces to a
 $d\times d$  matrix $q$ with $d$ the rank of $E$, and $\Pi_{\theta,
   \theta^\prime}$ reduces to $\pi_{\theta, \theta^\prime}=\frac{1}{2i\pi}
\int_{C_{\theta, \theta^\prime, r}} \lambda^{-1} \, q\,
(q-\lambda)^{-1} \, d\lambda$ which is a finite dimensional projector. Hence,  the generators 
 $\left[L_{\theta, \theta^\prime}^Q(t)\right]$ reduce to  generators $[l_{\theta,
   \theta^\prime}^q(t)]=[e^{2i\pi\, t \pi_{\theta, \theta^\prime}(q)}]$
 built from  projectors 
 $\pi_{\theta, \theta^\prime}(q)$.
\end{rk}
{\bf Proof:} We take the proof from \cite{LP} \footnote{As pointed out to us
  by R. Ponge, the proof can probably be shortened using results of \cite{BL}
  to show directly that $K_0(Cl^0(M, E))$ is generated by idempotents
  $\Pi_{\theta, \theta^\prime}(Q)$.}. For any algebra ${\mathcal A}$, let $K_0({\mathcal A})$ denote the group
of formal differences of homotopy classes of idempotents in the direct limit
${\rm gl}_\infty(A)$ of matrix algebras ${\rm gl}_n(A)$\footnote{A natural
  embedding ${\rm gl}_n(A)\to {\rm gl}_{n+1}(A)$ of  an $n\times n$
matrix $a\in {\rm gl}_{n}(A)$ in ${\rm gl}_{n+1}(A)$  is obtained  inserting $a$ in the upper left corner
and filling the last line and column with zeroes.}. 
When ${\mathcal A}$ is a {\it good topological} algebra \cite{Bo}, the Bott
periodicity isomorphism: 
\begin{equation}\label{bott}
  \begin{array}{ccc}
   K_0({\mathcal A}) &\longrightarrow & \pi_1\left({\rm GL}_\infty( {\mathcal A}\right)) \\
           \lbrack P \rbrack     & \longmapsto & e^{2i\pi t P} 
    \end{array}
\end{equation}
holds.  
Since for any vector bundle $E$ over $M$, the algebra $Cl^{0}(M, E)$
is a good topological algebra (which essentially boils down to the
fact that the inverse of a classical pseudodifferential operator  remains a
classical pseudodifferential operator), applying (\ref{bott}) to
${\mathcal A}=Cl^{0}(M, E)$ reduces the proof down  to checking  that 
$K_0(Cl^{0}(M, E))$ is generated by idempotents 
$\Pi_{\theta, \theta^\prime}(Q)$. \\
  The  exact sequence $$ 
   0\longrightarrow Cl^{-1}(M,E)\longrightarrow Cl^0(M,E) \stackrel{\sigma^L}{\longrightarrow} 
   C^\infty(S^*M,p^*(\End E)) \longrightarrow 0 
 $$ where $p: S^*M\to M$ is the canonical projection of the cotangent sphere
 to the base manifold $M$,
 gives rise to a long exact sequence in $K$-theory: 
\begin{equation}\label{longexactseq}
\begin{array}{ccccc} 
 K_0(Cl^{-1}(M,E))& \rightarrow & K_0(Cl^0(M,E)) &\stackrel{\sigma^L_0}{\rightarrow} &
    K_0(C^\infty(S^*M,p^*(\End E))) \\
 \uparrow {\rm Ind} & & & & \downarrow 0 \\
  K_1(C^\infty(S^*M,p^*(\End E))) &\stackrel{\sigma^L_1}{\leftarrow}& K_1(Cl^0(M,E))
  & \leftarrow & K_1(Cl^{-1}(M,E))=0.  
\end{array}
\end{equation}
On the other hand, on the grounds of results of Wodzicki, $K_0(C^\infty(S^*M,p^*(\End E)))$ is generated by the classes $\pi_{\theta,
  \theta^\prime}(\sigma_L(Q))$ where as before $\sigma_L(Q)$ is the leading
symbol of an elliptic operator $Q\in \Cl(M, E)$; this combined with the surjectivity of the map
$\sigma_0^L$ in the diagram (\ref{longexactseq}) yields the result. \endsquare
\\ \\ Higher homotopy groups were derived in \cite{BW} and \cite{R} from which
we quote some results without proofs \footnote{Even homotopy groups were also
  described in \cite{R} leading to further results which we do not report on
  here. Also, the statement we quote here holds provided one allows for bundles with
arbitrary large rank. }.
 \begin{prop} For odd $k$\begin{enumerate}
\item \cite{BW}(Proposition 15.4)
$$
\pi_k\left(\Cl_{Id}^{0, *}(M,E)\right)\simeq\quad  \Z,$$
where  $\Cl_{Id}^{0, *}(M,E):=\{A\in\Cl^{0, *}(M,E), \quad \sigma_L(A)=Id\},$ 
\item \cite{R} (Theorem 1)
$$
\pi_k\left(\Cl^{0, *}(M,E)\right)\simeq\quad  K_0(C^\infty(S^*M,p^*(\End E))),$$
which is therefore generated by the homotopy classes of the loops
given by (\ref{eq:generators}).
\end{enumerate}
\end{prop}

\section{A class of infinite dimensional  manifolds}

We consider a  class  of infinite dimensional manifolds and vector bundles
inspired from the geometric setup of index theory \cite{B}, \cite{BGV} and 
close to those introduced in 
\cite{P1} (under the name of weighted  manifolds and bundles) and further used
in \cite{CDMP}, \cite{PR1,PR2} (under the name of $\pdo$-manifolds and bundles). 
It consists of Fr\'echet vector bundles with typical fibre $\Ci(M, E)$ for
some reference finite rank  vector bundle $\pi: E\to M$ over $M$. \\ 
Consider a smooth fibration $\M\to X$ of smooth manifolds modelled
  on a closed manifold $M$ and a
  fibre bundle $\pi:\E\to \M$ over $\M$ with typical fibre $E\to M$.
\begin{rk} In
  the context of the family index theorem,  $\E= \F\otimes
  \vert\Lambda_\pi\vert$  for some vector bundle  $\F\to \M$ and $\Lambda_\pi
  $ is the vertical density bundle which, when restricted to the fibres of
  $\M$ may be identified with the bundle of densities along the fibre.
\end{rk}
 Let us denote by 
  $\pi_* \E\to X$ the infinite dimensional Fr\'echet bundle with fibre
  $\Ci(M_x, E_{\vert_{ M_x}})$ over
  $x\in X$ modelled on $\Ci(M, E)$ with $M$ the model fibre of  $ \M$ and $E$
  the model fibre of  $\E$.
\begin{defn}We call  a Fr\'echet vector bundle admissible if it is of the form  $\pi_* \E$ for some
  finite rank vector bundle $\E\to \M$ over a smooth fibration $\M\to X$ of smooth
  closed manifolds.\\ We call a  Fr\'echet manifold admissible if its
  tangent bundle is an  admissible vector bundle.
\end{defn}
 \begin{rk}Locally, over an open subset $U\subset X$, $$\E_{\vert_U}\simeq U\times
   M\times V$$ where $V$ is a finite vector space and $M$ a closed manifold. A
   change of local trivialisation of the fibration $\M_{\vert_U}\simeq U\times M$
   induces a diffeomorphism $f: M\to M $ in 
  $ {\cal D}(M)$ whereas a change of local trivialisation of the finite rank
  vector bundle $\E_{\vert_U\times M}\simeq U\times M\times V$  induces a
  transformation in ${\rm Gl}(V)$.
\end{rk}
\begin{rk} When $M$ reduces to a point $\{*\}$ then an admissible vector bundle
  reduces to a finite rank vector bundle over $X$ modelled on some vector
  space $V$ with transition maps in ${\rm Gl}(V)$.
\end{rk}
\begin{ex}\label{ex:mappingspaces} Let $N$ be a Riemannian manifold, then the space $X:=\Ci(M, N)$ of
  smooth maps from $M$ to $N$ is a
  Fr\'echet manifold with tangent space at point $\gamma$ given by $\Ci(M,
  \gamma^*TN)$. The tangent bundle $T\Ci(M, N)$ can therefore be realised as
  $\pi^*\E$ where $\M\to X$ is the trivial fibration with fibre $M$ and
  $\E$ the vector bundle over $X\times M$ with fibre at $(\gamma, m)\in \Ci(M,
  N)\times M $
  given by the vector space $$\E_{(\gamma, m)}= \gamma^*T_{\gamma(m)}N
  \quad{\rm so }\quad {\rm that} \quad \pi_* \E_{\gamma}=
 \Ci(M,  \gamma^*T_{\gamma}N).$$ Hence $\Ci(M, N)$ is an admissible manifold.\\
In passing note  that this  manifold, which is modelled on $\Ci(M, \R^n)$ where $n$
  is the dimension of $N$ can be equipped with an atlas induced by the
  exponential map  $\exp^N$ on $N$, a local chart being of the type 
$\phi_\gamma(u)(x)= \exp^N_{\gamma(x)}(u(x))$. The transition maps are
multiplication operators. 
\end{ex}
\begin{ex}\label{ex:mappingroups} In particular,  mapping groups $ \Ci(M, G)$ with $G$ a finite
  dimensional Lie group are admissible Fr\'echet manifolds. \\
The  left  action $l_g: y\mapsto g\cdot y$ on $G$ induces a left action $L_g:
\gamma\mapsto g\cdot \gamma$ on $\Ci(M, G)$ and a vector field $V(\gamma)\in
\Ci(M,\gamma^*T_\gamma G)$ is left-invariant if $\left(L_g\right)_*
V(\gamma)=V(g\cdot \gamma)$ for all $\gamma\in \Ci(M, G)$. Left invariant
vector fields on $\Ci(M, G)$ can be identified with  elements of the Lie
algebra $\Ci(M, {\rm Lie}(G))$.
\end{ex}
\begin{ex}\label{ex:diffeomorphisms} The group 
$${\cal D}(M):=\{f\in C^\infty(M, M), \quad \exists f^{-1}\in C^\infty(M,
M)\}$$ of smooth diffeomorphisms of $M$
is a Fr\'echet Lie group \cite{O} (see \cite{N} for a review) with Lie algebra $\Ci(M, TM)$ where $TM$ is the
tangent bundle to $M$. It is an admissible Fr\'echet manifold since its
tangent bundle $\bigcup_{f\in {\cal D}(M) }\Ci(M,  f^*TM)$ can be realised as a bundle
$\pi^*\E\to X$ with $X= {\cal D}(M)$ and where $\E$ is a vector bundle over the trivial
fibration $\M=X\times M$ with
fibre above $(f, M)$ given by  the bundle $f^*TM$. 
\end{ex} 
We also introduce a  class of connections inspired from the ones arising in
the family index geometric setup and similar to the ones considered in \cite{P1} later
named $\pdo$-connections  in
\cite{PR1,PR2}.  \\ Following \cite{Sc},
\cite{PS2}, let for any $a\in \C$,  $\Cl^a(\M,
\E)$ denote the bundle over $X$ with fibre over $x$ given by $\Cl^a\left(M_x,
\E_{\vert_{M_x}}\right)$ so that locally, above an open subset $U$ of $X$ we have 
$$\Cl^a(\M,
\E)_{\vert_U}\simeq U\times \Cl^a(M, E).$$ Let $\Cl(\M, \E)$ be the bundle of
algebras generated by $ \bigcup_{a\in
  \C}\Cl^a(\M, \E)$.
\begin{rk} When $M$ reduces to a point $\{*\}$ and $\E\to \M$ reduces to a
 finite rank vector  bundle $E\to X$,  then $Cl^a(\M, \E)= \Cl(\M,
  \E)={\rm End}(E)$ the endomorphism bundle over $M$.
\end{rk}
\begin{ex} For mapping spaces as described in Example \ref{ex:mappingspaces}
  the fibre of the bundle $\Cl(\M, \E)$  above $\gamma
$ is given by
$$\Cl(\M, \E)_{\gamma}= \Cl(M, \gamma^*TN).$$
\end{ex}
\begin{ex} For mapping groups as described in  Example \ref{ex:mappingroups},
we can specialise to left-invariant $\pdo$s, $A(\gamma)\in \Cl(M,
\gamma^*TG)$ such that $A(g\cdot \gamma)\circ{ L_g}_* ={ L_g}_* \circ
A(\gamma)$ for all $g\in \Ci(M, G)$ and $\gamma\in \Ci(M, G)$.
\end{ex}
Above an open
set $U\subset X$, a connection on a finite rank  vector bundle $E\to X$ is locally  of the form:
$$\nabla_V \sigma= d\sigma(V)+ \theta^U(V)\sigma, \quad\forall V\in T_xX,
\quad \forall \sigma\in \Ci(U, E)\quad {\rm with}\quad \theta^U(V)\in {\rm
  End}(E_x).$$ 
In view of the generalisation from ${\rm End}(E_x)$ to
$\Cl^0(M_x,\E_{\vert_{M_x}})$, we introduce a class of connections locally of
the form:
\begin{equation}\label{eq:pdoconnection} \nabla_V\sigma = d\sigma(\tilde V)+ \theta^U(\tilde V)\quad\forall V\in T_xX
\quad \forall \sigma\in \Ci(X, \pi_*\E) \quad {\rm with}\quad \theta^U(\tilde
V)\in \Cl^0(M_x, \E_{\vert_{M_x}})
\end{equation}
with $\E$ a finite rank bundle  over a fibration $\pi:\M\to X$
of manifolds equipped
with some  horizontal distribution $V\in T_xX \mapsto \tilde V\in T_{(x, m)}\M$.\\ 
Similar classes of connections  were  considered  in \cite{P1} and
later in \cite{PR1,PR2} under the name of $\pdo$-connection.
\begin{defn}
A connection $\nabla$  on an admissible Fr\'echet vector  bundle  $\pi_*\E \to X$ with
$\E$ a finite rank bundle  over a fibration $\pi:\M\to X$
of manifolds is admissible whenever
\begin{enumerate}
\item   the connection $\nabla$ is induced by a connection on a finite rank
  vector bundle:
\begin{equation}\label{eq:connectiononfibration}\nabla_V= \nabla^{\E}_{\tilde V}
\end{equation}
for some connection $\nabla^\E$ on the finite rank vector bundle $\E\to \M$
and some  horizontal   distribution $V\in
  T_xX \mapsto \tilde V\in T_{(x, m)}\M$  on
  $\M$,

\item or when the fibration $\M= X\times M$ is trivial, if  locally over an open
  subset $U\subset X$ 
$$ \nabla = d+ \theta^U \quad {\rm with}\quad \theta^U\in \Omega^1(U,  \Cl^0(\M, \E)).$$
\end{enumerate}

\end{defn}
\begin{rk}The two conditions have a  non void intersection; indeed, in the case of a  horizontal distribution on a trivial fibration
 $\M\simeq X\times M$  the first
  condition locally reads $\nabla= d+\theta^U$ with $\theta^U$ given by a
  multiplication operator valued one form on $U$.  
\end{rk}
\begin{rk}

Admissible connections fulfill condition (\ref{eq:pdoconnection}); we  choose
the trivial horizontal lift in the case of a trivial connection.
\end{rk}
\begin{rk} When
$M$ reduces to a point $\{*\}$ then any connection is admissible 
since $\E\to \M$ boils down to a bundle $E\to X$ and locally $\nabla=
d+\theta^U$ with 
$\theta^U\in\Omega^1(U, {\rm End}( E))$.
\end{rk}
\begin{lem} An admissible  connection on an admissible bundle $\pi_*\E\to X$ has
  curvature in $\Omega(X, \Cl^1(\M, \E))$. When the fibration $\M\to X$ is
  trivial with trivial distribution then it lies in  $\Omega(X, \Cl^0(\M, \E))$.
\end{lem}
{\bf Proof:}
\begin{enumerate}
\item If locally,  $\nabla= d+\theta^U$ with $\theta^U\in \Omega^1(U, \Cl^0(\M,
  \E))$  then $\Omega= d\theta^U+ \theta^U\wedge\theta^U$ lies in  $\Omega^2(U, \Cl(\M,
  \E))$. Since the curvature is a globally defined two form,  $\Omega$ lies in
  $\Omega^2(X, \Cl^0(\M, \E))$.
\item If $\nabla_V= \nabla^\E_{\tilde V}$ then \begin{eqnarray*}
\Omega(U, V)&=&[ \nabla_U,  \nabla_V
]- \nabla_{[U, V]}\\
&=&[ \nabla^\E_{\tilde U},  \nabla^\E_{\tilde V}
]- \nabla^\E_{\widetilde{[U, V]}}\\
&=&  \Omega^{ \E}(\tilde U, \tilde
V)-\nabla^\E_{\widetilde{[U, V]}-[\tilde U, \tilde V]}\\
&=&  \Omega^{ \E}(\tilde U, \tilde
V)-\nabla^\E_{T(U, V)}\\
\end{eqnarray*}
where  $T(U, V)= \widetilde{[U, V]}-[\tilde U, \tilde V]$ is the
curvature of the connection on $\M$. Since $\Omega^{ \E}(\tilde U, \tilde
V)$ is a multiplication operator, it follows that $\Omega(U, V)$ is a first
order differential operator which therefore lies in $\Cl^1(\M, \E)$.  When the
distribution on $\M$  is trivial, then $T(U, V)=0$ and  $\Omega(U, V)$
lies in $\Cl^0(\M, \E)$.
\end{enumerate}
\endsquare

\begin{ex} \label{ex:connectiononcurrentgroups}On mapping groups considered in
  Example \ref{ex:mappingroups}, we restrict to 
  left-invariant connections i.e. connections $\nabla$  such that if $V, W$ are
  left-invariant then so is $\nabla_VW$ left-invariant. An admissible
  left-invariant connection is defined by a left-invariant one-form
 $
\theta_0\in  \Omega^1(X, \Cl^0(M,{\rm Lie}(G))).$    \\ With the notations
  of Example \ref{ex:mappingroups}, let $(M, g)$ be a closed
  Riemannian manifold and let  $G$ be a semi-simple Lie group of compact type \footnote{This ensures
  that the Killing form is non degenerate and that the adjoint representation
  ad on the Lie algebra is antisymmetric for this bilinear form}.  
  \\ Let $Q_0:=
  \Delta\otimes 1_{{\rm Lie}(G)}$ where $\Delta$ stands for the Laplace
  Beltrami operator on $M$.  D. Freed in \cite{F} introduces a family of left-invariant one-forms  parametrised by $s\in \R$ on the $H^s$-Sobolev closure $H^s(M, G)$ of
  the the mapping group $\Ci(M, G)$ (see formula (1.9) in \cite{F}):
\begin{eqnarray*}
\theta_0^s(V)&:= &\frac{1}{2} \left({\rm ad}_V + (Q_0+ \pi_0)^{-s} {\rm ad}_V
(Q_0+ \pi_0)^{s}\right. \\
& -& \left. (Q_0+ \pi_0)^{-s}{\rm ad}_{ (Q_0+ \pi_0)^{-s}V}\right)\quad
\forall V\in \Ci(M, {\rm Lie}(G)).
\end{eqnarray*}
Here $\pi_0$ stands for the orthogonal projection onto the kernel of $Q_0$
which is finite dimensional.\\
These give rise to  left-invariant connections $\nabla^s$ on
$H^s(M, G)$ 
 which in turn induce 
connections\footnote{They are weak since they are defined by
  weak metrics on $L^2(M, {\rm Lie}(G))$; they are not determined by
  the usual six term formula \cite{F}.} on the mapping group  $\Ci(M, G)$. Since $\theta_0^s\in \Omega^1(
\Ci(M, G), \Cl^0(M, {\rm Lie}(G))$, these define $\pdo$-connections;
only if $s=0$ does this one-form correspond to a multiplication operator. \\
The curvature $\Omega^s$ is given by a left-invariant two-form:
\begin{equation}\label{eq:curvaturecurrentgroup}
\Omega_0^s(U, V)= [\theta_0^s(U), \theta_0^s(V)]-\theta_0^s([U V])
\end{equation} and by \cite{F} (Proposition 1.14), the map
$$R^s(U, V): W\mapsto \Omega^s(W, U) \, V$$
is a pseudodifferential operator of order max$(-1, -2s)$. It lies in
$\Omega^2\left(M, \Cl^0(M,{\rm Lie}(G))\right)$.
 \end{ex}
\begin{rk}
 These results were extended
 to loop spaces in \cite{MRT}, where it was shown that the curvature
 operator $R^s$ on $\Ci(S^1, N)$ built in a similar manner has order at most
 $-1$ for $s>\frac{3}{2}$. 
\end{rk}
\begin{ex}\label{ex:connectiononfibrations}
Connections of the type (\ref{eq:connectiononfibration}) arise in the
geometric setup underlying the index
theorem for families \cite{B} and in  determinant bundles
associated with families of Dirac operators \cite{BF}.  Such connections are in fact
slightly perturbed  by adding  the divergence of the
horizontal lift w.r. to the Riemannian volume element on the fibres of $\M\to
X$ in order to produce
unitary connections:
$$\tilde \nabla^{\E}_V= \nabla^{\pi_*\E}_{\tilde V}+ {\rm div}_M(\tilde V).$$
\end{ex}
\section{Singular Chern-Weil forms in infinite dimensions}
 We aim at generalising Chern-Weil formalism to admissible vector bundles,
 defining when possible traces  $\lambda(\nabla^{2i})$ of even powers of
 admissible connections. In view of the two types of admissible connections  we  distinguish  two cases:
\begin{enumerate}
\item  the curvature lies in  $\Omega^2\left(X, \Cl^0(\M, \E)\right)$,
\item  the curvature  lies in  $\Omega^2\left(X, \Cl^1(\M, \E)\right)$. 
\end{enumerate}
Accordingly, to build traces of powers of the curvature along the lines of the first section, by the results of
section 3  we  have at our disposal two types of (singular) traces, namely
traces that vanish on smoothing operators:
\begin{enumerate}
\item leading symbol traces  (\ref{eq:leadingsymboltrace}) and the
  noncommutative residue  (\ref{eq:Wodzickiresidue})  on $\Cl^0(M, E)$
\item the noncommutative residue (\ref{eq:Wodzickiresidue}) on $\Cl(M, E)$.
\end{enumerate}

Let  $V$ be  the model space of $\E\to \M$. Since a change of trivialisation
$((x,m), v)\in \M\times V\mapsto ((x, m), C\, v)\in \M\times V, \quad C\in {\rm
  Gl}(V)$  of $\E\to \M$ induces a change of trivialisation 
$$(x,A)\in X\times \Cl(M, V)\mapsto (x, C^{-1} \, A\, C)\in X\times \Cl(M,
V),\quad C\in \Ci(M, {\rm Gl}(V))$$ of the bundle  
$\Cl(\M, \E)\to X$, we need to make sure the traces we implement are invariant
under the adjoint action of $\Ci(M, {\rm Gl}(V))$. The following lemma gives
more, namely their invariance under the action of invertible zero order $\pdo$s. 
\begin{lem}Let $E\to M$ be a finite rank vector bundle over a closed manifold
  $M$. For any $ C\in \Cl^{0,*}(M, E)$ \begin{equation}\label{eq:covariantres}{\rm res}(C^{-1}\,A\,
C)={\rm res}(A), \quad \forall A\in \Cl(M,
E)\end{equation} 
and \begin{equation}\label{eq:covarianttrleadingsymb}{\rm tr}^\tau(C^{-1}\,A\,
C)={\rm tr}_0^\tau(A)\quad  \forall A\in \Cl^0(M,
E).\end{equation} 
\end{lem}
\begin{rk} Here, as pointed out in the introduction, the group $\Cl^{0,*}(M, E)$ of invertible classical $\pdo$s generalises  the structure group ${\rm GL}(V)$.
\end{rk}
{\bf Proof:} Both properties follow from the cyclicity of the respective
traces.\endsquare\\ \\
Just as the ordinary trace on matrices induces a trace on the endomorphism
bundle  End$(E)$ of any finite rank vector bundle $E\to X$, the noncommutative residue
(resp. and the leading symbol traces) therefore induce a noncommutative residue
(resp. and  leading
symbol traces) on bundles $\Cl(\M, \E)$ (resp. on  $\Cl^0(\M, \E)$)
associated to any finite rank vector bundle $\E\to \M$ over a fibration of
manifolds $\M\to X$. \\
As in the finite dimensional case,  to a form $ \alpha(x)= A(x)\, dx_1\wedge \cdots\wedge dx_d$ in
$\Omega\left(X,\Cl(M,
E)\right)$ (resp.
$\Omega\left(X,\Cl^0(M,
E)\right)$)    corresponds  a form
  ${\rm res}(\alpha)(x):={\rm res}(A(x))\, dx_1\wedge \cdots\wedge dx_d$
(resp. and ${\rm tr}_0^\tau(\alpha)(x):={\rm tr}_0^\tau(A(x))\, dx_1\wedge
\cdots\wedge dx_d$  for any $\tau\in \Ci(S^*M)^\prime$) in
$\Omega(X)$. 
As in the finite dimensional case (see (\ref{eq:dtr}) and (\ref{eq:coboundarytr})),  we first check that  the linear form
$\lambda= {\rm res}$ (resp.  $\lambda={\rm
  tr}_0^\tau$ for any $\tau\in \Ci(S^*M)^\prime$)  on $\Omega(X, {\cal A})$ with ${\cal A}= \Cl(\M,\E)$  (resp.   ${\cal A}=\Cl^0(\M,\E)$) obeys the following properties 
\begin{equation}\label{eq:dlambda} \left[d,
    \lambda\right](\alpha):=d\,\lambda (\alpha)- \lambda (d\,
\alpha)=0\quad  \forall \alpha\in \Omega\left(X,{\cal A}\right)\end{equation}
and 
\begin{equation}\label{eq:coboundarylambda} 
\partial \lambda(\alpha, \beta):= \lambda\left(\alpha \wedge
  \beta+(-1)^{\vert \alpha\vert \, \vert \beta\vert} \beta \wedge
  \alpha\right)=0 \quad \forall \alpha, \beta\in \Omega\left(X, {\cal A}\right).
\end{equation} The first property is  easily   checked from the very
definition of the two types of traces which involve integrals over the
cotangent unit sphere of the manifold $M$ of the trace of a homogeneous part
of the symbol of the operator. The second property is a direct consequence of
their cyclicity. \\
As in the finite dimensional case, we  then  infer the following  lemma.

\begin{lem}\label{lem:nablalambda}For any $\alpha\in \Omega\left(X,\Cl(\M, \E)\right)$
\begin{equation}\label{eq:nablares}[\nabla, {\rm res}](\alpha):= d\,  {\rm
    res}(\alpha)-{\rm res}([\nabla,\alpha])=0.
\end{equation}
For any $\alpha\in \Omega\left(X,\Cl^0(\M, \E)\right)$, for any $\tau\in \Ci(S^*M)^\prime$
\begin{equation}\label{eq:nablatr0}[\nabla, {\rm tr}_0^\tau](\alpha):= d\,  {\rm
    tr}_0^\tau(\alpha)-{\rm tr}_0^\tau([\nabla,\alpha])=0.
\end{equation}
\end{lem}
{\bf Proof:} The proof goes as for Lemma \ref{lem:nablatr} using properties
(\ref{eq:dlambda}) and (\ref{eq:coboundarylambda}).\endsquare\\ The
constructions carried out in section 1 then go through leading to the
following result. 
\begin{thm} Let ${\cal E}=\pi^* \E \to X$ be an admissible vector bundle with
  $\E\to \M$ a finite rank vector bundle over a fibration $\M\to X$, equipped
  with a connection with curvature $\Omega$. For any $i\in \N$, the $i$-th residue
  Chern-Weil form
\begin{enumerate}
\item  
$ {\rm res}\left(\Omega^{i}\right)$   is closed with de Rham cohomology class independent of the choice of
connection,
\item if $\Omega \in \Omega^2(X, \Cl^0(\M, \E))$ then 
$ {\rm tr}_0^\tau \left(\Omega^i\right)$   is also closed with de Rham cohomology class independent of the choice of
connection.
\end{enumerate}
We call these {\bf singular Chern-Weil classes}. 
\end{thm}
{\bf Proof:}
The proof goes as in Proposition \ref{prop:ChernWeilfinitedim} using 
Lemma \ref{lem:nablalambda}.\endsquare
\begin{rk} Singular Chern-Weil forms are clearly insensitive to smoothing
  perturbations of the connection.
\end{rk}
\begin{ex} We refer to \cite{RT} where  an explicit example of an
infinite rank bundle   with non vanishing first residue Chern class is built,
i.e. where the noncommutative residue is used as an Ersatz for the usual trace
on matrices to build a first Chern form with non vanishing de Rham class; it is a
bundle over the two-dimensional sphere  $S^2$ with fiber modelled on Sobolev sections of a trivial
line bundle over the three dimensional torus $T^3$.\end{ex}
\begin{ex}Going back to Example \ref{ex:connectiononfibrations}, we saw that 
$\Omega^{\pi_* \E}$ is a differential operator valued two form so that its noncommutative
residue vanishes. When the fibration  is trivial  $\M\simeq M\times X$, $\Omega^{\pi_* \E}$ reduces to a
multiplication operator so that its leading symbol traces are well-defined and
the two forms 
$$(U, V)\mapsto {\rm tr}_0^\tau\left( \widetilde \Omega^{\pi^*\E}(U,V)\right)=
{\rm tr}_0^\tau\left(  \Omega^{\E}( \tilde U, \tilde V)\right)
$$
can give rise to non trivial singular Chern-Weil classes.
\end{ex}
But unfortunately, singular Chern-Weil classes generally  seem too coarse to capture
interesting information since  most   examples   lead to vanishing singular
Chern-Weil classes.

\begin{ex}\label{ex:resChernformcurrentgroup} Going back to Example \ref{ex:connectiononcurrentgroups}, it was shown in \cite{F} that $W\mapsto \Omega_0^s(W, U)V$ is
  conditionally trace-class, i.e. that ${\rm tr}_{{\rm Lie}(G)}\Omega_0^s(\cdot,
  U)V$ is trace-class, from which we infer that the residue vanishes. \\
For the same reason, the leading symbol traces also vanish.\\ Thus, singular
Chern-Weil forms vanish in the case of mapping groups. 
\end{ex}

\begin{rk} In \cite{MRT}, the authors  actually used the fact that  
 residue  Pontrjagin forms  vanish on loop spaces $\Ci(S^1, N)$  as the starting point to build
  singular Chern-Simons classes. They focus on manifolds $N$ with stably trivial
  tangent bundle. 
\end{rk}

\section{Weighted Chern-Weil forms; discrepancies}
 We now address the issue of how to  build analogs of Chern-Weil forms using 
extensions of the ordinary trace of the type
discussed in section 4 instead of singular traces on $\pdo$s.
 We show how  one stumbles on various discrepancies inherent to the
 fact that these extensions do not actually define traces; however, it is
 instructive to  describe the obstructions to carrying out these constructions
 in order to later find ways to circumvent them. \\ \\
We first recall a result of Wodzicki, Guillemin and popularised by Kontsevich
and Vishik in  \cite{KV} which relates  the complex  residue of the
canonical trace of a holomorphic family at a pole with the noncommutative residue of the
family at this
pole: for any $C\in \Cl(M, E)$ and any holomorphic family $C(z)\in \Cl(M, E)$ with order $\alpha(z)$  such that $C(0)=C$ and
$\alpha^\prime(0)\neq 0$, the following identity holds: 
\begin{equation}\label{eq:KVopclassical}
 {\rm
  Res}_{z=0} {\rm TR}\left(C(z)\right)= -\frac{1}{ \alpha^\prime(0)}{ \rm
res}(C).
\end{equation}
This formula provides ways to measure various defects of
weighted traces.\\
\subsection{The Hochschild coboundary of a weighted trace}
 The first defect  is an obstruction to its cyclicity. The Hochschild
 coboundary of a
 linear form $\lambda$  on $\Cl(M, E)$ is defined by $$\partial \lambda(A, B):= \lambda\left([A, B]\right).$$ The following proposition says this coboundary is
 local (see e.g. \cite{CDMP},\cite{Mi1,Mi2},  \cite{MN}). 
\begin{prop}  Let $Q\in \Cl(M, E)$ be an admissible operator of positive order $q$
  with spectral cut $\theta$.
Let $A\in \Cl (M, E)$, $B\in (M, E)$
then 
\begin{equation}\label{eq:Hochschild}\partial {\rm tr}_\theta^Q\left(A, B\right)=-\frac{1}{q}\, {\rm res}\left( A\,
  [B,\log_\theta Q]\right).
\end{equation}
\end{prop}
{\bf Proof:} Using the vanishing of the canonical trace on non integer order
brackets we can write
\begin{eqnarray*}
{\rm TR}\left([A\,, B]\,Q_\theta^{-z}\right)&=&{\rm TR}\left(A\, B\,
  Q_\theta^{-z}- B\, A\, Q_\theta^{-z} \right)\\
&=&{\rm TR}\left(A\, B\,
  Q_\theta^{-z}- A\, Q_\theta^{-z}\,B\right)\\
&=& {\rm TR}\left(A\,[ B,
  Q_\theta^{-z} ] \right).
\end{eqnarray*}
 The family $C(z):=\frac{A\,[ B,
  Q_\theta^{-z} ]}{z}\in \Cl(M, E)$ is a holomorphic family of order $a-b-
q\, z$ and $C(0)=- A\,[ B,
  \log Q ]$.  By (\ref{eq:KVopclassical}) we get:
\begin{eqnarray*}
{\rm tr}_\theta^Q\left([A, B]\right)&=& {\rm fp}_{z=0}{\rm TR}\left(A\, B\, Q_\theta^{-z}- B\, A\, Q_\theta^{-z}\right)\\
&=& {\rm fp}_{z=0}{\rm TR}\left(A\,[ B,
  Q_\theta^{-z} ]\right)\\
&=&  {\rm Res}_{z=0}{\rm TR}\frac{\left(A\,[ B,
  Q_\theta^{-z} ]\right)}{z}\\
&=& -\frac{1}{q}{\rm res}\left(A\, [B, \log_\theta Q]\right).
\end{eqnarray*}
\endsquare
\subsection{Dependence on the weight}
 Formula (\ref{eq:KVopclassical}) also provides  a way to measure the
dependence on the weight $Q$. We first need a technical lemma.
\begin{lem} Let $Q\in \Cl(M,E)$ be admissible of order $q>0$ and spectral cut
  $\alpha$ and let  $A\in \Cl(M, E)$. Then
$${\rm tr}_\theta^{Q^t}(A)= {\rm tr}_\theta^Q(A)\quad \forall t>0.$$ 
\end{lem}
{\bf Proof:} We write ${\rm TR}(A\, Q_\theta^{-z})= \frac{a_{-1}}{z}+
a_0+o(z)$
in which case we have ${\rm TR}\left(A\, \left(Q^t_\theta\right)^{-z}\right)= {\rm
  TR}\left(A\, \left(Q_\theta\right)^{-t\, z}\right) =\frac{a_{-1}}{tz}+
a_0+o(t\,z)$ so that $${\rm tr}_\theta^{Q^t}(A)={\rm fp}_{z=0}{\rm TR}(A\,
\left(Q^t_\theta\right)^{-z})=a_0={\rm tr}_\theta^{Q}(A).$$
\endsquare\\ \\
The following proposition provides   a well-known expression of the dependence
on the weight \cite{KV},\cite{Ok}.
\begin{prop}\label{prop:dependenceQ}Let $Q_1, Q_2\in \Cl(M, E)$ be two admissible operator with
  positive orders $q_1, q_2$
  and spectral cuts $\theta_1, \theta_2$.
Let $A\in \Cl (M, E)$, then
$${\rm tr}_{\theta_1}^{Q_1}(A)-{\rm tr}_{\theta_2}^{Q_2}(A)= {\rm res}\left( A\,
  \left(\frac{\log_{\theta_1} Q_2}{q_2}-\frac{\log_{\theta_2} Q_1}{q_1}\right)\right).$$
\end{prop} 
{\bf Proof:} For simplicity, we leave out the explicit mention of the spectral
cut. Applying formula (\ref{eq:KVopclassical}) to the family
$C(z):= \frac{A}{z}\,\left(
    Q_1^{-\frac{z}{q_1}}- 
    Q_2^{-\frac{z}{q_2}}\right)$ which is a holomorphic family of classical
  operators of   order $a-z$ with $C(0)= A\, \left(\frac{\log Q_2 }{q_2}-
  \frac{\log Q_1}{q_1}\right)$ we write
\begin{eqnarray*}
{\rm tr}^{Q_1}(A)-{\rm tr}^{Q_2}(A)&=&{\rm tr}^{Q_1^{\frac{1}{q_1}}}(A)-{\rm tr}^{Q_2^{\frac{1}{q_2}}}(A)\\
&=& {\rm Res}_{z=0} {\rm TR}
\left( \frac{A\,\left(Q_1^{-\frac{z}{q_1}}- Q_2^{-\frac{z}{q_2}}\right)}{z}\right)\\
&=& {\rm res}\left( A\, \left( \frac{\log Q_2}{q_2}- \frac{\log Q_1}{q_1}\right)\right).
\end{eqnarray*}
\endsquare
\subsection{Exterior differential of a weighted trace}

\begin{prop}\label{prop:diffweightedtrace}\cite{CDMP}, \cite{P1} Let $Q\in C^\infty(X,\Cl(M, E))$ be a
  differentiable family  (for the topology described previously) of operators of
  fixed order $q$
  and spectral cut $\theta$
  parametrised by a manifold $X$.
Let $A\in \Cl (M, E)$, then the trace defect $\left[d,{\rm tr}^{Q}\right]:=d\, {\rm tr}^{Q}-{\rm
  tr}^{Q}\circ d$ is local as a noncommutative residue:
\begin{equation}\label{eq:exteriordiff}\left[d,{\rm tr}_\theta^{Q}\right](A)= -\frac{1}{q}\, {\rm res}\left( A\,
 d\, \log_{\theta} Q\right).
\end{equation}
\end{prop} 
{\bf Proof:} Again, for convenience, we drop the explicit mention of the spectral cut. Let $h\in \Ci(X,TX)$ be a smooth vector field then by 
Proposition \ref{prop:dependenceQ} we have
\begin{eqnarray*}
d{\rm tr}^{Q}(A)(h)&=&\lim_{t\to 0}\frac{ {\rm tr}^{Q+dQ(th)}(A)-{\rm tr}^{Q}(A)}{t}\\
&=& -\frac{1}{q}\,\lim_{t\to 0}\frac{ {\rm res}\left( A\,
  \left(\log (Q+dQ(th))-\log Q\right)\right)}{t} \\
&=& -\frac{1}{q}\,\frac{ {\rm res}\left(\lim_{t\to 0}\left( A\,
  \left(\log (Q+dQ(th))-\log Q\right)\right)\right)}{t} \\
&=& -\frac{1}{q}\, {\rm res}\left( A\,
  d \,\log Q(h)\right), \\
\end{eqnarray*}
where we have used the continuity of the noncommutative residue on operators
of order $a={\rm ord } A$ since $ A\,
  \left(\frac{\log (Q+dQ(th))-\log Q}{t}\right)$ has order $a$ for any $t>0$.
\endsquare

\subsection{Weighted traces extended to admissible fibre bundles}
The covariance property  (\ref{eq:covarianttrace}) generalises to weighted
traces as follows. 
\begin{lem} \cite{P1}, \cite{PS2}  Order, ellipticity,  admissibility and spectral cuts are preserved under
  the adjoint action. For any admissible operator $Q\in \Cl(M, E)$ with
  spectral cut $\alpha$, for any operator
  $A\in \Cl(M, E)$ we have \begin{equation}\label{eq:covariance}{\rm
      tr}_\theta^{{\rm ad}_CQ }({\rm ad}_CA)= {\rm tr}_\theta^Q(A)\quad\forall C\in \Cl^{
  *}(M, E).\end{equation}
\end{lem}
{\bf Proof:} For simplicity, we drop the explicit mention of the spectral cut.\\Since the leading symbol is multiplicative we have 
$$\sigma_L({\rm ad}_CQ)= {\rm ad}_{\sigma_L(C)}\sigma_L(Q)$$ from which it
follows 
 that order, ellipticity, admissibility and spectral cuts are preserved by the
 adjoint action.\\
Let us  observe that 
\begin{eqnarray*}
\left( {\rm ad}_CQ\right))^{-z}&=&\frac{1}{2i\pi}\int_\Gamma \lambda^{-z} (\lambda- {\rm
  ad}_CQ)^{-1} \\
&=& \frac{1}{2i\pi}\int_\Gamma \lambda^{-z} {\rm
  ad}_C(\lambda- {\rm
  ad}_CQ)^{-1} \\
&=& {\rm ad}_CQ^{-z}
\end{eqnarray*}
Consequently,  
\begin{eqnarray*}
{\rm tr}^{C^{-1} Q C}(C^{-1}A C)&=&{\rm fp}_{z=0} {\rm TR}\left(C^{-1}A C\,
\left(C^{-1}Q C\right)^{-z}\right)\\
&=& {\rm fp}_{z=0} {\rm TR}(C^{-1}A C\,
C^{-1}Q^{-z} C)\\
&=& {\rm fp}_{z=0} {\rm TR}(C^{-1}A Q^{-z} C)\\
&=& {\rm fp}_{z=0} {\rm TR}(A Q^{-z} )\\
&=& {\rm tr}^Q(A)
\end{eqnarray*}
where we have used the fact that the canonical trace vanishes on non integer
order brackets.
\endsquare\\\\
Since the adjoint action ${\rm ad}_C: A\mapsto C^{-1} AC$  of $\Cl^{0, *}(M,
E)$ on $\Cl(M, E)$ preserves the spectrum and the invertibility of the leading symbol, it   makes sense to define
the subbundle  $ \Ell^{\rm adm}(\M, \E))$ of $\Cl(\M, \E)$ of fibrewise
admissible elliptic $\pdo$s with spectral cut $\theta$; since it also preserves the order we can define
$\Q$ to be a smooth admissible elliptic section of order $q$ of $\Cl(\M,
\E)$ in which case $\Q(x)\in \Cl(M_x, \E_{\vert_{M_x}})$ and we set:
$${\rm tr}_\theta^\Q(A)(x):= {\rm tr}_\theta^{\Q(x)}(A(x) )\quad \forall A\in
\Ci(X, \Cl(\M, \E)),\quad \forall x\in X.$$
$\Q$-weighted traces can further be extended  to  forms $ \alpha(x)= A(x)\, dx_1\wedge \cdots\wedge dx_d$ in
$\Omega\left(X,\Cl(\M,
\E)\right)$  by ${\rm tr}_\theta^\Q(\alpha)(x):={\rm tr}_\theta^\Q(A(x))\, dx_1\wedge
\cdots\wedge dx_d$ and using linearity.
\subsection{Obstructions to closedness of weighted Chern-Weil forms} 
\begin{thm}\label{thm:nablaweightedtrace} \cite{CDMP} An admissible  connection $\nabla$
  on an admissible vector bundle $\pi_* \E$ induces a connection $[\nabla, A]:= \nabla\circ
  A-A\circ \nabla$ on $\Cl(\M, \E)$. For any $\alpha\in
  \Omega^p\left(X, \Cl(\M, \E)\right)$ and any admissible section $\Q\in \Ci\left(X, {\cal A}\right) 
$ with constant spectral cut $\theta$ and constant order $q>0$, the trace defect $[\nabla, {\rm tr}_\theta^Q]:=d\, {\rm tr}_\theta^\Q(\alpha)- {\rm
  tr}_\theta^\Q\left([\nabla, \alpha]\right) $ is local and explicitely given by:
\begin{equation}\label{eq:nablaweightedtrace}[\nabla, {\rm tr}_\theta^\Q](\alpha)= \frac{(-1)^p}{q} {\rm res} \left(
  \alpha\,
[\nabla, \log_\theta \Q]\right)\quad \forall \alpha\in \Omega^p(X,
\Cl(\M, \E)).
\end{equation}
\end{thm}
{\bf Proof:} We prove the result for a zero form and drop the explicit mention
of the spectral cut for simplicity. The result easily extends to
higher order forms. In a local trivialisation of ${\cal E}$ over an open subset $U$ of $X$ we write
$\nabla= d+\theta$ so that $[\nabla, \cdot]= d+[\theta, \cdot]$.
In this local trivialisation we have for any $\in \Ci(X, \Cl(\M,\E))$:
\begin{eqnarray*}
[\nabla, {\rm tr}^Q](A)&=& d\,\left( {\rm tr}^Q(A)\right)- {\rm
  tr}^Q\left(\left[ \nabla, A\right]\right)\\
&=&  d\, \left({\rm tr}^Q(A)\right)- {\rm
  tr}^Q\left(d\,A\right)\, - {\rm
  tr}^Q\left(\left[ \theta, A\right]\right)\\
&=& -\frac{1}{q} {\rm res}\left( A\,d \log Q\right) -\frac{1}{q} {\rm
  res}\left( A\,[\theta,\log Q]\right)\\
&=& -\frac{1}{q} {\rm res}\left( A\,[\nabla , \log Q]\right)
\end{eqnarray*}
where we have combined (\ref{eq:Hochschild}) and (\ref{eq:exteriordiff}).

\endsquare
\\ \\The obstruction $[\nabla^{\rm ad}, {\rm tr}^Q]$ described in  Theorem
\ref{thm:nablaweightedtrace} prevents  a straightforward
generalisation  of the Chern-Weil formalism to an infinite dimensional setup where the trace on
matrices is replaced by a weighted trace provided   the connection are
admissible  connections.
\begin{cor} \label{cor:nonclosed}Let $\nabla$ be an admissible connection on ${\cal E}=\pi_*\E\to
  X$ with curvature $\Omega$ and let $\Q$ be an admissible section of $\Cl(\M, \E)$ with spectral cut
  $\theta$ and  constant positive order $q$. Then 
$$d\, {\rm tr}_\theta^\Q(\Omega^i)= \frac{1}{q} {\rm res} \left( \Omega^i\, 
[\nabla, \log_\theta \Q] \right).$$
\end{cor}
{\bf Proof:} We follow the finite dimensional  proof (see Proposition \ref{prop:ChernWeilfinitedim}).
\begin{eqnarray*}
 d\, {\rm tr}_\theta^\Q(\Omega^i)&=& [ d\, {\rm tr}_\theta^\Q](\Omega^i)+ {\rm
   tr}_\theta^\Q([\nabla, \Omega^i ])\\
&=& [ d\, {\rm tr}_\theta^\Q](\Omega^i)\quad {\rm by} \quad {\rm the}\quad {\rm Bianchi} \quad {\rm identity}\\
&=& \frac{1}{q} {\rm res} \left( \Omega^i\, 
[\nabla, \log_\theta \Q] \right)\quad{\rm by}\quad
(\ref{eq:nablaweightedtrace}).
\end{eqnarray*}
\endsquare
\section{Renormalised Chern-Weil forms on $\pdo$ Grassmannians}
In view of Corollary \ref{cor:nonclosed} which tells us that a weighted
trace of a power of the curvature is generally not closed, it  seems hopeless
to use weighted 
traces as a substitute for ordinary traces in order to extend finite
dimensional Chern-Weil formalism to infinite dimensions. However, there are
different ways to circumvent this difficulty, one of which is  to introduce
counterterms in order to compensate for the lack of  closedness measured in
Corollary \ref{cor:nonclosed} by a noncommutative
residue. Such a renormalisation procedure by the introduction of counterterms
can be carried out in a hamiltonian approach to  gauge theory as it was shown
in joint work with J. Mickelsson \cite{MP} on which we report here.\\
Let us   first review a finite dimensional situation which will serve as a
model for infinite dimensional generalisations.\\ We consider  the finite-dimensional Grassmann manifold ${\rm Gr}(n,n)$
consisting of rank $n$ projections in
$\C^{2n}$, which we   parametrise by grading
operators $F= 2P -1,$ where $P$ is a finite rank projection.
\begin{lem} The even forms 
\begin{equation}\label{eq:finitedimomegaj} \omega_{2j} = \tr\left( \, F(dF)^{2j}\right), \end{equation}
where $j=1,2,\dots$ are closed  forms on ${\rm Gr}(n,n)$.
 \end{lem}
\begin{rk}The cohomology of  $G_n(\C^{2n}):={\rm Gr}(n, n)$ is 
known to actually be generated by even (nonnormalised) forms of the type
$\omega_{2j}, j=1, \cdots,
n$.  This follows from the fact (see Proposition 23.2 in \cite{BT}) that the Chern classes of the quotient bundle
$Q$ over $G_n(\C^n)$ defined by the exact sequence
$0\to S\to G_n(\C^{2n})\times \C^n\to Q\to 0$ where $S$ is the universal bundle over $G_n$ with
fibre $V$ above $V$, generate the cohomology ring
$H^*(G_n(\C^n))$.  Then the $j$-th Chern class of $Q$ turns out to be
proportional to ${\rm tr}(F(dF)^{2j})$ where $P(V)$
stands for the orthogonal projection on $V$. 
\end{rk}
{\bf Proof:} By the traciality of $\tr$ we have
\begin{eqnarray}
 d\, \omega_{2j} &=& d\, \tr\left( \, F(dF)^{2j}\right)\nonumber\\
&=& \tr\left( (dF)^{2j+1}\right)\nonumber\\
&=& \tr\left( F^2\, (dF)^{2j+1}\right)\nonumber\\
&{}&{\rm since} \quad F^2=1\nonumber\\
&=& -\tr\left( F\, (dF)^{2j+1}\, F\right)\nonumber\\
&{}&{\rm since} \quad F\, dF=-dF \, F\nonumber\\
&=& -\tr\left( (dF)^{2j+1}\, F^2\right)\nonumber\\
&{}&{\rm since} \quad \tr([A, B])=0\nonumber\\
&=& -\tr\left( (dF)^{2j+1}\right)\nonumber\\
&=& 0.
\end{eqnarray}
\endsquare
\\ \\
We now want to extend these constructions to  $\pdo$ Grassmannians.\\  
Let us  consider a finite rank  bundle $ \E$ over a trivial fibration
 $\pi: \M=M\times X \to X$ with typical fibre a closed (Riemannian) spin
 manifold $M$. Let $D_x\in \Cl(\M, \E), x\in X$ be a smooth family of Dirac operators
 parametrised by $X$.
\\ On each open subset $U_\lambda:=\{x \in x,\lambda\notin {\rm
  spec}(D_x)\}\,  \subset X$ there is a well defined map
\begin{eqnarray*}
F: X&\to & \Cl^0(M, E)\\
x &\mapsto &
F_x:=(D_x-\lambda I)/|D_x -\lambda I|.
\end{eqnarray*} 
 Since $F_x^2= F_b$, $P_x:= \frac{I+ F_x}{2}$ is a projection, the
  range Gr$(M, E):= {\rm Im} F$ of $F$ coincides with   the Grassmannian consisting of classical pseudodifferential  projections $P$ with kernel and
cokernel of infinite rank, acting in the complex Hilbert space $H:= L²(M,
E)$ of  square-integrable sections of the vector bundle $E$ over the
compact manifold $M.$\\
  Since the map $x\mapsto F_x$  is generally not contractible  we want
to define from $F_x$  cohomology classes on $X$ in the way we built Chern-Weil
classes in finite dimensions.  This issue usually arises 
in hamiltonian quantization in field theory, when  the physical 
space $M$ is an odd dimensional manifold. In this infinite dimensional setup  traces are  generally 
ill-defined, so that we use weighted traces as in the previous section. As
expected, there are a priori obstructions to the closedness of the corresponding
weighted forms.
\begin{prop} Let $Q\in \Cl(M, E)$ be a {\bf fixed} admissible elliptic operator with
  positive order. The exterior differential of the form 
\begin{equation} \label{eq:trQomegaj}
\omega^Q_{2j}(F) = \tr^{Q} \left(
    F(dF)^{2j}\right)
 \end{equation} 
on  ${\rm Gr}(M, E)$:
$$  d\omega^Q_{2j}= \frac{1}{{2q}} {\rm res} \left([\log Q,F] (dF)^{2k+1}
  F\right).$$ 
 is a local expression which  only  depends 
on $F$ modulo smoothing operators. 
 \end{prop}
{\bf Proof:} The locality and the dependence on $F$ modulo smoothing operators
 follow from the expression of the exterior differential in terms of a
Wodzicki residue.  To derive this expression, we mimic the finite dimensional proof, taking into account that
this time ${\rm tr}^Q$ is not cyclic:
\begin{eqnarray*}
  d\omega^Q_{2j} & = & d{\rm tr}^{Q}\left( F(dF)^{2j}\right) \\
&=& \tr^Q\left( (dF)^{2j+1}\right)\nonumber\\
&=& {\rm tr}^{Q}\left( F^2 (dF)^{2j+1}\right)\\
&=& -{\rm tr}^{Q} \left(F (dF)^{2j+1} F \right)\\
&{}& {\rm since} \quad F\, dF=-dF \, F\\
&=&\frac{1}{q} {\rm res}\left( [\log Q, F] (dF)^{2j+1}F\right) -
{\rm tr}^{Q}\left( (dF)^{2j+1} F^2\right) \\ 
&=&\frac{1}{q} {\rm res}\left( [\log Q, F] (dF)^{2j+1}F\right) -
{\rm tr}^{Q}\left( (dF)^{2j+1}\right), 
 \end{eqnarray*}
where we have used (\ref{eq:Hochschild}) to write 
$$ {\rm tr}^{Q} \left([F, (dF)^{2j+1} F ] \right)= -\frac{1}{q} {\rm res} \left(F\,[ (dF)^{2j+1} F, \log
  Q]\right)=\frac{1}{q} {\rm res} \left([  F, \log
  Q]\, (dF)^{2j+1}F\right).$$
Hence 
$${\rm tr}^{Q}\left( F^2 (dF)^{2j+1}\right)= \frac{1}{2q} {\rm res}\left(
  [\log Q,
F] (dF)^{2j+1}F\right)$$ from which  the result then  follows. 
\endsquare  \\ \\ 
 Let us consider the  map \begin{eqnarray*}
\sigma:  X &\to &  \Cl^0(M, E)/ Cl^{-\infty}(M, E) \\
x&\mapsto & \bar F(x):= p\circ F(x)
\end{eqnarray*}
 where $p:\Cl^0(M, E)\to \Cl^0(M, E)/ Cl^{-\infty}(M, E)$ is the canonical projection
 map. 
\\ \\
The following theorem builds  from the original forms
$\omega^Q_{2j}$ new ``renormalised'' forms which are closed in contrast with the original ones.
\begin{thm}\label{thm:renormgauge} When $\sigma(X)$ is contractible, there are  even forms $\theta^Q_{2j}$ such that 
$$ \omega^{{ren}, Q}_{2j} := \omega^Q_{2j} - \theta^Q_{2j} $$ 
is closed. 
The forms $\theta^Q_{2j}$ vanish when  the order of $(dF)^{2j+1}$ is less than
-dim$\,M$. This holds in particular if  the order of $(dF)^{2j}$ is less than -dim$\,M$ in which
case  $ \omega^{ren, Q}_{2j} = \omega^Q_{2j}= {\rm tr}(F(dF)^{2j})$ is independent of $Q$.
\end{thm}
{\bf Proof:} The form  $d\omega^Q_{2j}$ being a Wodzicki residue, it is
insensitive to smoothing perturbations  and is therefore  a pull-back  by the
  projection map $p$  of a form
  $\beta_{2j}^Q$. 
The  pull-back of $\beta_{2j}^Q$  with respect to $\sigma$ 
is a closed form $\theta^Q_{2j+1}$ on $X$ which is exact since $\sigma$ is
contractible. Indeed, selecting a contraction $\sigma_t$ with $\sigma_1 =\sigma$ and $\sigma_0$ a constant map, 
we have the standard formula $d\theta^Q_{2j} =
\theta^Q_{2j+1},$ with \begin{equation} \theta_{2j} = \frac{1}{2j+1} \int_0^1 t^{2j}\iota_{\dot\sigma_t}  \theta^Q_{2j+1}(\sigma_t) dt. \end{equation} where $\iota_X$ is the contraction by a vector field $X$ and the dot means differentiation with respect 
to the parameter $t$.\\
 When the order of $(dF)^{2j+1}$ is less than -dim$\,M$ the correction terms $\theta_{2j}^Q$ 
vanish and if  the order of $(dF)^{2j}$ is less than -dim$\,M$, the weighted trace ${\rm tr}^Q$ coincides with the usual trace so that
the naive expression $\omega^Q_{2j}$ is a closed form independent
of $Q$.
\endsquare
\\ \\ This way, on builds renormalised Chern classes $[\omega^{{\rm ren},
  Q}_{2j}]$. We refer the reader to \cite{MP} for
the two form case which arises in the
quantum field theory gerbe \cite{CMM}.
\section{Regular Chern-Weil forms  in infinite dimensions}

We describe further geometric setups for which weighted traces actually do
give rise to closed Chern-Weil forms.  \\ \\
 Mapping groups studied by Freed \cite{F} and
later further investigated in e.g. \cite{CDMP}, \cite{M}, \cite{MRT} provide a
first illustration of such a situation. \\ 
Going back to Example \ref{ex:mappingroups}, we specialise to the
circle $M=S^1$; the Sobolev based loop group $H_e^{\frac{1}{2}}(S^1, G)$  can be equipped
with a complex structure and its first Chern form  was studied by Freed
\cite{F}. We saw in Example  \ref{ex:resChernformcurrentgroup} that the
corresponding curvature is conditionally trace-class which leads to the
following result. 
\begin{prop}\cite{F} Theorem 2.20  (see also \cite{CDMP} Proposition 3)
Let $Q_0\in \Cl(S^1)$ be an admissible elliptic operator on $S^1$ with spectral
cut $\theta$. Let $\nabla^{\frac{1}{2}}$ be a left-invariant connection on
$\Ci(S^1, G)$ with curvature given by a two form $\Omega_0^{\frac{1}{2}}$ as in (\ref{eq:curvaturecurrentgroup}),
 then by  (\ref{eq:condtraces}) we have 
$${\rm tr}_\theta^{Q_0}(\Omega_0^{\frac{1}{2}})={\rm tr}_{\rm
    cond}(\Omega_0^{\frac{1}{2}})= {\rm TR}(\Omega_0^{\frac{1}{2}}). $$ It defines a closed form which coincides with Freed's
    conditioned first Chern form. 
\end{prop}
\begin{rk}It was observed by Freed in \cite{F} that this weighted first Chern form ${\rm
    tr}_\theta^{Q_0}(\Omega_0^{\frac{1}{2}})$  relates to the K\"ahler form on
  the based loop group $H_e^{\frac{1}{2}}(S^1, G)$. See also \cite{CDMP} for
  further interpretations of this two-form. 
\end{rk}
 Another  way around the  obstructions described previously is to choose a weight $\Q$ and a
connection $\nabla$ such that the bracket $[\nabla, \log_\theta \Q] $
vanishes; this can be achieved using superconnections, leading to a second
geometric setup in which regularised traces do give rise to closed Chern-Weil
forms. 
\begin{defn}A superconnection (introduced by Quillen  \cite{Q}, see also  \cite{B},
 \cite{BGV}) on an admissible vector bundle $\pi_* \E$ where $\pi:\M\to X$ is
 a fibration of manifolds, adapted to a smooth
 family of
formally self-adjoint elliptic $\pdo$s $D\in C^\infty\left(X, \Cl^d(\M, \E)\right)$  with odd
parity is a linear map $\A$ acting on  $\Omega
\left(X,
\pi_* \E\right)$ of odd parity with respect
 to the
$\Z_2$-grading such that:
$$\A(\omega\cdot \sigma)= d\omega \wedge \sigma +
(-1)^{\vert \omega\vert}\omega\wedge  \A(\sigma)\quad \forall
\omega \in \Omega(X), \sigma \in \Omega\left(X,
\pi_* \E \right) $$
and  $$\A_{[0]}:= D,$$
where we have written  $\A= \sum_{i=0}^{{\rm dim} B}
\A_{[i]}$ and
$\A_{[i]}:  \Omega^* \left(X, {\mathcal
E}\right)\mapsto \Omega^{*+i} \left(X, {\mathcal
E}\right).$
\end{defn}
\begin{ex} An admissible connection $\nabla$ as in
  (\ref{eq:connectiononfibration}) gives rise to  a superconnection
$$\A:= \nabla+ D.$$
\end{ex}
The curvature  of a superconnection $\A$ is  a $\pdo$-valued form 
 $\A^2\in \Omega^2 \left(X, \Cl(\M, \E)\right)$; it actually is a differential
 operator valued two form. Since
$\A^2= D^2+\A^2_{[>0]}$ where $ \A^2_{[>0]}$ is a $\pdo$-valued form of positive
degree, just as $D^2$,  $\A^2$ is elliptic and admissible.  \\
Following \cite{Sc} and \cite{PS2}, we call a $\pdo$-valued form  $\omega= \sum_{i=0}^{{\rm dim} B}
\omega_{[i]}$ with $\omega_{[i]}:  \Omega^* \left(X, \pi_* \E\right)\mapsto \Omega^{*+i} \left(X, \pi_* \E\right)$ elliptic, resp. admissible, resp. with spectral cut $\alpha$
whenever $\omega_{[0]}\in \Cl(X,\pi_*\E)$ has these properties. We refer the reader to \cite{PS2} for detailed
explanations on this point.\\
 Since  $\A^2_{[0]}= D^2$, 
$\A^2$ is an  elliptic $\pdo$-valued form with spectral cut $\pi$, and hence 
an
admissible $\pdo$-valued form. Its complex powers and logarithm can be defined
as for ordinary  admissible  $\pdo$s. \\
With these conventions, weighted traces associated with  fibrations of $\pdo$-algebras can be
generalised to weights given by  admissible $\pdo$-valued forms such as the
curvature $\A^2$ of the superconnection. Along the lines of the proof of the
previous theorem one can check that the trace defect $[\A, {\rm tr}^{\A^2}]$ vanishes:
 \begin{prop}\label{prop:superconnection} \cite{PS2} For any $\omega\in\Omega\left(X, \Cl\left({\cal E}\right)\right)$ we have
$$d\, {\rm str}^{\A^2}(\omega)={\rm str}^{\A^2}([\A, \omega]).$$
\end{prop}
{\bf Proof:} Since
$$d\, {\rm str}^{\A^2}(\omega)={\rm str}^{\A^2}([\A, \omega])+ [\A,{\rm
  str}^{\A^2}] (\omega) $$
this  follows from Theorem \ref{thm:nablaweightedtrace} combined
with the fact that $[\A, \log \A^2]=0$. \endsquare
\begin{thm} The form
${\rm str}^{\A^2}(\A^{2j})$ defines a closed form  called  the $j$-th Chern form  
associated with the superconnection $\A$ and a de Rham cohomology
class independent of the choice of
connection.
\end{thm}
{\bf Proof:}  This follows from  the Bianchi identity $[\A,\A^{2j}]=0$
combined with Proposition \ref{prop:superconnection}.\endsquare 
\begin{rk} Let $\A= D+ \nabla$ be a
superconnection associated with a family of Dirac operators $D$ on a trivial
fibration of manifolds. It was observed in \cite{MP}   that the expression
$${\rm tr}^{D^2} (\nabla^{2j})- {\rm tr}^{\A^2}(\A^{2j})_{[2j]}$$ 
-which compares the naive infinite dimensional analog  ${\rm tr}^{D^2}
(\nabla^{2j})$ of the finite dimensional Chern form ${\rm tr} (\nabla^{2j})$
and the 
 closed form ${\rm tr}^{\A^2}(\A^{2j})_{[2j]}$ built from the super
 connection-  is local in as far as it  is  insensitive to smoothing
 perturbations of the connection. The weighted Chern-Weil form ${\rm
   tr}^{\A^2}(\A^{2j})_{[2j]}$ is therefore
 interpreted as a renormalised version of ${\rm tr}^{D^2}
(\nabla^{2j})$.   This is similar to  the formula derived in the previous
section  where   a residue correction
term was added to the naive weighted form involving the curvature.
\end{rk}
If we specialise to  a fibration $\pi: \M \to X$ of
even-dimensional closed spin manifolds with   the  Bismut superconnection $\A:=
D+\nabla +c(T)$,  $c$
 the Clifford multiplication and $T$ the curvature of the horizontal
distribution on $\M$, we get an explicit description of the $j$-th Chern
form  ${\rm str}^{\A^2} (\A^{2j})$.  Indeed, as a consequence of  the  local index theorem for
families \cite{B} (see also \cite{BGV}),   the component of degree $2j$ of the form 
${\rm str}^{\A^2}(\A^{2j})$ can be expressed in terms of the  $\hat A$-genus $\hat A(\M/X)$  on the vertical fibre 
  of $\M\to X$ and the Chern character  ${\rm ch}(\E_{\M/
X})$ on the
  restriction of $\E$ to the vertical fibre.
\begin{thm} \cite{MP} Let $\A$ be a superconnection adapted to a family of
  Dirac operators on even dimensional spin manifolds parametrised by $X$, then 
\begin{equation}\label{eq:{eq:strChernj}} {\rm str}^{\A^2} \left( \A^{2j}\right)_{[2j]}= \frac{(-1)^j j!}{(2i\pi)^{\frac{n}{2}}} \, \left( \int_{\M/X} \hat A(\M/X) \wedge
  {\rm ch}(\E_{\M/B}) \right)_{[2j]}.\end{equation}
\end{thm}
{\bf Proof:} Along the lines of the heat-kernel proof of the index theorem
(see e.g. \cite{BGV})  we introduce the kernel $k_\e ( \A^{2})$ 
 of $e^{-\e \A^2}$ for some $\e>0$. Since $D$ is a family of Dirac
   operators, we have   (see e.g. chap. 10 in \cite{BGV}) 
\begin{equation} \label{eq:heatkernelexp}
k_\e(\A^{2})(x, x)\sim_{\e \to 0} \frac{1}{(4\pi \e)^\frac{n}{2}} \sum_{j=0}^\infty \e^j k_j( \A^{2})(x,
x).
\end{equation} We  observe that the $j$-th Chern form associated with $\A$ is given by an integration along
  fiber of $\M:$ 
 $${\rm str}^{\A^2} \left( \A^{2j}\right)=  \frac{(-1)^j j!}{(4\pi )^{\frac{n}{2}}} \int_{\M/B}{\rm str}( k_{j+\frac{n}{2}}(
      \A^{2}))$$ and proceeed to compute ${\rm str}( k_{j+\frac{n}{2}}(
      \A^{2}))$.\\  Let us introduce  Getzler's rescaling which  transforms  a homogeneous form $\alpha_{[i]}$ of degree $i$ to the expression
$$\delta_t \cdot \alpha_{[i]}\cdot \delta_t^{-1}=  \frac{\alpha_{[i]}}{ \sqrt t^i},$$ so that a superconnection
$\A=\A_{[0]}+\A_{[1]}+\A_{[2]} $ transforms to 
$$\tilde \A_t= \delta_t \cdot \A\cdot \delta_t^{-1}=\A_{[0]} +
\frac{\A_{[1]}}{\sqrt t}+\frac{\A_{[2]}}{t}. $$
As in   \cite{BGV} par. 10.4, in view of 
the asymptotic expansion (\ref{eq:heatkernelexp})
we have:
\begin{eqnarray*}
{\rm ch}(\A_t)&=& \delta_t \left( {\rm str} (e^{- t \A^2})\right)\\
&\sim_{t \to 0}& (4 \pi t)^{-\frac{n}{2}} \sum_j t^j\int_{\M/B} \delta_t \left({\rm str}( k_j (\A^2))\right)\\
&\sim_{t \to 0} & (4 \pi )^{-\frac{n}{2}} \sum_{j, p} t^{j-(n+p)/2}\left( \int_{\M/B}{\rm str} \left(   k_j (\A^2)\right)\right)_{[p]},\\
\end{eqnarray*} 
so that 
\begin{equation}\label{eq:rescaledkernelasymp}
{\rm fp}_{t=0}{\rm ch}(\A_t)_{[p]}= (4 \pi )^{-\frac{n}{2}}  \left(
  \int_{\M/B}{\rm str} \left(   k_{\frac{p+n}{2}} (\A^2)\right)\right)_{[p]}.
\end{equation}
The family index  theorem  \cite{B} (see also  Theorem 10.23 in \cite{BGV})
yields the existence of the  limit as $t\to 0$ and 
$$\lim_{t\to  0} {\rm ch}(\A_t)= (2 i\pi)^{-\frac{n}{2}} \, \int_{\M/B} \hat A (\M/B)\wedge
  {\rm ch}(\E_{\M/B}).$$
Combining these two facts leads to:
$$\left( \int_{\M/B}{\rm str} \left(   k_{\frac{n+2j}{2}} (\A^2)\right)\right)_{[2j]}=\frac{(4 \pi )^{\frac{n}{2}}}{(2i\pi)^{\frac{n}{2}}} \left( \int_{\M/B} \hat A(\M/B) \wedge
  {\rm ch}(\E_{\M/B}) \right)_{[2j]}.$$
It follows that
 $$ {\rm str}^{\A^2} \left( \A^{2j}\right)_{[2j]}=  \frac{(-1)^j j!}{(2i\pi)^{\frac{n}{2}}} \, \left(\int_{\M/B} \hat A(\M/B) \wedge
  {\rm ch}(\E_{\M/B})\right)_{[2j]}.$$ 
\endsquare\\ \\
As could be  expected,  $\A^2$-weighted Chern forms therefore relate to the Chern character $(2i\pi)^{\frac{n}{2}} \,  \int_{\M/B} \hat A(\M/B) \wedge
  {\rm ch}(\E_{\M/B})$ of a family
of Dirac operators associated with the fibration $\M\to X$ \cite{BGV}.
\vfill \eject 



\printindex
\end{document}